# IMMERSED BOUNDARY METHOD FOR THE COMPLETE ELECTRODE MODEL IN ELECTRICAL IMPEDANCE TOMOGRAPHY

JÉRÉMI DARDÉ, NIAMI NASR, AND LISL WEYNANS

ABSTRACT. We propose an immersed boundary scheme for the numerical resolution of the Complete Electrode Model in Electrical Impedance Tomography, that we use as a main ingredient in the resolution of inverse problems in medical imaging. Such method allows to use a Cartesian mesh without accurate discretization of the boundary, which is useful in situations where the boundary is complicated and/or changing. We prove the convergence of our method, and illustrate its efficiency with two dimensional direct and inverse problems.

## 1. INTRODUCTION

This work is dedicated to the numerical resolution of direct and inverse problems related to Electrical Impedance Tomography (EIT) using an immersed boundary method. The aim of EIT is to reconstruct the electrical conductivity distribution inside a domain by imposing electrical currents on the boundary of this domain, and measuring the resulting voltages on the same boundary. It has several applications in the medical field, in particular in lung monitoring and stroke detection [1, 8, 30]. Mathematically, the problem, known as *Calderón problem* or *inverse conductivity problem*, is a severely ill-posed inverse problem. We refer to [5, 9, 49] and the references therein for an overview on the inverse conductivity problem.

In practical experiments, the currents are driven inside the body of interest through a collection of surface electrodes, no current being driven between the electrodes. For each current pattern, the potential differences between the electrodes are measured. This practical setting is accurately modeled by the Complete Electrode Model (CEM) [11, 48]. It takes into account the shape of the electrodes as well as the shunting effect, that is the thin resistive layer that appears at the interface between the electrodes and the object during the measurements. The CEM is known to correctly predict experimental data, and therefore is widely used in the numerical resolution of both direct and inverse problems related to EIT.

Usually, Finite Element Methods (FEM) are used to compute solutions of the CEM. This leads to fast and efficient numerical methods in simple geometries. However, when the geometry is complicated, the FEM can become very expensive, especially in three-dimensional problems. This is due to the need of an accurate discretization of the geometry, in particular on the electrodes where the potential varies rapidly in space, imposing the use of very refined meshes.

Another difficulty could arise when the geometry is not perfectly known, which usually leads to highly incorrect reconstructions [22, 42]. Iterative Newton-type algorithms that reconstruct both the conductivity and the geometry have been proposed to tackle this difficulty [15, 16]. In such algorithms, the geometry changes from one iteration to the other. Therefore, a new mesh is created at each step, leading again to a costly numerical implementation.

In the present paper, our objective is to propose an alternative to the numerical resolution of CEM with classical mesh-adapted numerical methods, by using *Immersed Boundary Methods* (IBM), with the aim of reducing the cost of computation in complex and moving geometries. In such approach, the domain of interest is included in a larger domain with a simple geometry, typically a square in two-dimensional settings, or a cube in three-dimensional settings. The geometry of the domain of interest is encoded as the zero level-set of a known function, and is not discretized exactly. Instead, the larger domain is discretized using typically a Cartesian mesh. Of course, this approach requires a particular attention to include the boundary conditions of the initial problem.

The first method on Cartesian grids for elliptic problems was designed by Mayo in 1984 [39], and developed further in [40] and [38]. In that work an integral equation was derived to solve elliptic interface problems with piecewise coefficients to second order accuracy in maximum norm. Then LeVeque and Li (1994) [36] devised the very well known Immersed Interface Method (IIM), that relies on Taylor expansions of the solution on each side of the interface, with a local coordinate transformation near the interface to express the jump conditions in an appropriate frame. The elliptic operator is discretized on each grid point near the interface with formulas accounting for the jumps across the interface. In order to find these formulas, a linear system with six unknowns needs to be solved for each of the concerned grid points. Numerous developments of the IIM have been performed thereafter, among them [3, 37, 52]. Another class of Cartesian method was introduced by Zhou et al. is the Matched Interface and Boundary (MIB) method [53]. The solution on each side of the interface is extended on fictitious points on the other side. These fictitious values are computed by iteratively enforcing the lowest order interface jump conditions. Other classes of Cartesian methods also exist, less accurate in the case of interface problems, but probably simpler to implement: Gibou et al. [23, 24] developed methods inspired by Fedkiw's Ghost-Fluid Method [19] for multiphase flows. These methods are second order accurate for Dirichlet boundary conditions on arbitrary domains, but only first order accurate for interface problems. Then new methods were proposed to increase the accuracy in case of interface problems [4, 18].

In this paper, we approximate the CEM with a finite-difference scheme on a Cartesian mesh, and we follow the methodology of [13], where additional unknowns were defined at the intersections of the interface with the grid. These additional unknowns were used in the discretization of the elliptic operator near the interface, which avoided deriving specific finite differences formulas containing jump terms, or corrective terms. In order to solve the interface unknowns, the flux jump conditions were discretized and added to the linear system to solve. This use of additionnal unknowns at the interface between two domains, applied in the context of electrical impedance tomography, simplifies the discretization because the values at the boundary of the domain are directly involved in the equations modeling the input currents on the electrodes.

As stated above, the main difficulties for the discretization lie in the limit conditions. In particular, the limit conditions related to the input currents are integrals on the electrodes, which are non-standard from a finite differences point of view, and therefore need particular meticulousness. Another difficulty appears at the boundary of the electrodes, where the limit condition changes abruptly, leading to singularities in the solution. We tackle this difficulty by using the *smoothened Complete Electrode Model* [28], in which the transition between the conditions on electrodes and the condition in between electrodes in smooth, leading to smooth solutions.

Finally, as the problem is typically of *Neumann type* in order to accurately model the physical measurements setting, one has to deal with a *compatibility condition on the data* which is not easy to verify at the numerical level. We circumvent this last problem with a slight modification of the model, which allows to get rid of the compatibility condition without changing the solution of interest: this might be of interest on its own, as example to simplify numerical implementation of the CEM model by finite elements.

1.1. **Outline of the paper.** The paper is organized as follows. In Section 2, we present the Complete Electrode Model and its variations of interest for our study. In Section 3, we present the numerical scheme, and prove a linear rate of convergence. We present numerical results in a two-dimensional setting in Section 4, for both the direct problem of Electrical Impedance Tomography, the inverse conductivity problem, and the electrodes location problem. Some needed and technical computations are presented in the Appendix A.

## 2. The Complete Electrode Model for Electrical Impedance Tomography

2.1. **Geometrical and functional settings.** Here and in the following, we suppose that all functions and vectors are real-valued. Everything can be straightforwardly extended to the complex case.

Let $\Omega$ be a bounded domain of $\mathbb{R}^d$, $d = 2$ or $3$, with a smooth boundary $\partial\Omega$. We denote by $\nu$ the exterior unit normal of $\partial\Omega$. Let $(E_m)_{m=1,\ldots,M}$ be $M \geq 2$ mutually disjoint connected subdomains of $\partial\Omega$ with positive Lebesgue-measure. We set $E = \bigcup_{m=1}^M E_m$, and $E_c = \partial\Omega \setminus \overline{E}$.

We denote $\mathbb{R}^M_\diamond$ the subspace of mean-free vectors of $\mathbb{R}^M$,

$$\mathbb{R}^M_\diamond = \left\{ V \in \mathbb{R}^M,\ \sum_{m=1}^M V_m = 0 \right\},$$

and $\mathbb{R}^M_+$ the subset of $\mathbb{R}^M$ of vectors with positive components,

$$\mathbb{R}^M_+ = \left\{ V \in \mathbb{R}^M,\ V_m > 0 \quad \forall m \in \{1,\ldots,M\} \right\}.$$

For $s \in \mathbb{R}$, we introduce the quotient space $\mathscr{H}^s = (H^s(\Omega) \times \mathbb{R}^M)/\mathbb{R}$. Two elements $(u, U)$ and $(v, V)$ of $H^s(\Omega) \times \mathbb{R}^M$ are identified in $\mathscr{H}^s$ if they differ by an additive constant, in the following sense:

(2.1) $$(u, U) \sim (v, V) \Leftrightarrow \exists c \in \mathbb{R},\ u = v + c,\ U = V + c[1,\ldots,1].$$

Classicaly, $\mathscr{H}^s$ endowed with the norm

$$\|u, U\|_{\mathscr{H}^s} = \inf_{c \in \mathbb{R}} (\|u - c\|_{H^s(\Omega)} + \|U - c[1,\ldots,1]\|_{\mathbb{R}^M}),$$

is a Banach space.

For all $m$ in $\{1, \ldots, M\}$, let $e_m \subset E_m$ be of positive Lebesgue measure. We set

$$b : ((u, U), (v, V)) \in (H^1(\Omega) \times \mathbb{R}^M)^2 \mapsto \int_\Omega \nabla u \cdot \nabla v\, dx + \sum_{m=1}^M \int_{e_m} (u - U_m)(v - V_m)\, ds(x),$$

and, for all $(v, V) \in H^1(\Omega) \times \mathbb{R}^M$,

$$\|(v, V)\|_b^2 = b\left((v, V), (v, V)\right).$$

**Lemma 2.1.** *The bilinear form $b$ defines a scalar product on $\mathscr{H}^1$, with corresponding norm $\|.\|_b$ equivalent to the norm $\|.\|_{\mathscr{H}^1}$.*

*As a consequence, $(\mathscr{H}^1, b)$ is a Hilbert space.*

*Proof.* This result is standard, but we recall its proof for the reader's convenience. Firstly, we note that if $(u_1, U_1) \sim (u_2, U_2)$ and $(v_1, V_1) \sim (v_2, V_2)$ in the sense of (2.1), then $b\left((u_1, U_1), (v_1, V_1)\right) = b\left((u_2, U_2), (v_2, V_2)\right)$, so the Lemma makes sense.

Clearly $b$ is bilinear, symmetric and positive in $\mathscr{H}^1 \times \mathscr{H}^1$. Let $(v, V) \in \mathscr{H}^1$ be such that $\|v, V\|_b = 0$. Then $\nabla v = 0$ in $\Omega$, and there exists $c \in \mathbb{R}$ such that $v = c$ a.e. in $\Omega$. Furthermore

$$\int_{e_m} (v - V_m)^2 ds(x) = 0,\ \forall m \in \{1, \ldots, M\},$$

we obtain that $V = c[1, \ldots, 1]$, for the same constant $c$. Therefore, $(v, V) \sim (0, 0)$, and $b$ is a scalar product on $\mathscr{H}^1$.

Using the continuity of the trace application from $H^1(\Omega)$ to $L^2(\partial\Omega)$, it is readily seen that $\|.\|_b \lesssim \|.\|_{\mathscr{H}^1}$.

Conversely, let $(v_n, V_n)$ be a sequence in $\mathscr{H}^1$ such that $\|v_n, V_n\|_{\mathscr{H}^1} = 1$ and $\|v_n, V_n\|_b$ goes to zero as $n$ goes to infinity. Necessarily, there exists a sequence $(\tilde{v}_n, \tilde{V}_n)$ in $H^1(\Omega) \times \mathbb{R}^M$ such that

- $(v_n, V_n) \sim (\tilde{v}_n, \tilde{V}_n)$ in the sense of (2.1), and therefore $\|\tilde{v}_n, \tilde{V}_n\|_b$ goes to zero,
- $1 \leq \|\tilde{v}_n\|_{H^1(\Omega)}^2 + \|\tilde{V}_n\|_{\mathbb{R}^M}^2 \leq 2$, for all $n \in \mathbb{N}$.

As $(\tilde{v}_n)$ is a bounded sequence in $H^1(\Omega)$ and $(\tilde{V}_n)$ is a bounded sequence in $\mathbb{R}^M$, there exists a subsequence (still denoted $(\tilde{v}_n, \tilde{V}_n)$) and an element $(\tilde{v}_\infty, \tilde{V}_\infty) \in H^1(\Omega) \times \mathbb{R}^M$ such that $\tilde{v}_n$ weakly converges to $\tilde{v}_\infty$ in $H^1(\Omega)$ and $\tilde{V}_n$ strongly converges to $\tilde{V}_\infty$ in $\mathbb{R}^M$. As $\nabla \tilde{v}_n$ strongly converges to 0, $\tilde{v}_n$ actually strongly converges to $\tilde{v}_\infty$ in $H^1(\Omega)$, $\nabla \tilde{v}_\infty = 0$ and there exists a constant $c_\infty$ such that $\tilde{v}_\infty = c_\infty$ a.e. in $\Omega$.

Now, as for all $m \in \{1, \ldots, M\}$,

$$0 \leq \int_{e_m} (\tilde{v}_n - \tilde{V}_{n,m})^2 \, ds(x) \leq \|\tilde{v}_n, \tilde{V}_n\|_b^2,$$

taking the limit as $n$ goes to infinity leads to $\tilde{V}_{\infty,m} = \tilde{v}_{\infty|e_m} = c_\infty$, or equivalently $\tilde{V}_\infty = c_\infty[1, \ldots, 1]$. But as $(v_n, V_n) \sim (\tilde{v}_n, \tilde{V}_n)$, we have

$$1 = \|v_n, V_n\|_{\mathscr{H}^1} = \|\tilde{v}_n, \tilde{V}_n\|_{\mathscr{H}^1} \leq \|\tilde{v}_n - c_\infty\|_{H^1(\Omega)} + \|\tilde{V}_n - c_\infty[1, \ldots, 1]\|_{\mathbb{R}^M},$$

leading to a contradiction. Therefore $\|.\|_{\mathscr{H}^1} \lesssim \|.\|_b$, and the result follows. $\square$

2.2. **Generalized Complete Electrode Model.** We are now in position to introduce the generalized Complete Electrode Model we will use in the following. This model differs slightly from the standard Complete Electrode Model appearing in the context of Electrical Impedance Tomography, as we allow interior and boundary source terms, and we also allow spatially varying surface admittivities on the electrodes. Both features will be useful in the following study.

In our geometrical setting, $\Omega$ is the domain of interest, and each open set $E_m$ represents one of the electrodes. We denote $\sigma \in L^\infty(\Omega)$ the conductivity of the domain, which is assumed to verify the standard ellipticity condition

$$\sigma \geq c > 0 \text{ a.e. in } \Omega.$$

In the standard Complete Electrode Model, the thin resistive layer appearing at each electrode-body interface is represented by the vector of contact impedances $z \in \mathbb{R}_+^M$: one contact impedance $z_m$ per electrode $E_m$. The corresponding contact admittivity on electrode $E_m$ is then constant, equal to $\frac{1}{z_m}$. In our generalized model, we allow spatially varying admittivities. To do so, we consider $\xi_m \in \mathscr{Z}_m$, with

$$\mathscr{Z}_m = \{\xi \in L^\infty(E_m), \; \xi \geq 0, \; \xi \not\equiv 0\}.$$

These varying admittivities have been introduced in [28] as a way to obtain regularity on the solutions of the Complete Electrode Model. The standard Complete Electrode Model corresponds to the choice $\xi_m = \frac{1}{z_m}$, which is indeed in $\mathscr{Z}_m$.

Let $I \in \mathbb{R}^M$ represent the input currents imposed on the electrodes, and $(f, g) \in L^2(\Omega) \times L^2(\partial\Omega)$. The boundary value problem corresponding to the Complete Electrode Model reads: find $(u, U) \in \mathscr{H}^1$ such that

(2.2) $$\begin{cases} -\nabla \cdot (\sigma \nabla u) = f & \text{in } \Omega, \\ \sigma \nabla u \cdot \nu = g & \text{on } E_c, \\ \sigma \nabla u \cdot \nu + \xi_m(u - U_m) = g & \text{on } E_m, \; m = 1, \ldots, M, \\ \int_{E_m} \sigma \nabla u \cdot \nu \, ds(x) = I_m, & m = 1, \ldots, M. \end{cases}$$

Before proving that problem (2.2) is well-posed, we make several comments:

- A simple application of the divergence theorem shows that for problem (2.2) to admit a solution, it is necessary that $I$, $f$ and $g$ satisfy the *compatibility condition*

(2.3) $$\sum_{m=1}^M I_m + \int_\Omega f \, dx + \int_{E_c} g \, ds(x) = 0.$$

Of course, this is simply the well-known current conservation law applied to our system of equations. In the case of practical Electrical Impedance Tomography, where the source terms

$f$ and $g$ are null, we retrieve the mean-free condition on the input currents

$$\sum_{m=1}^{M} I_m = 0 \Leftrightarrow I \in \mathbb{R}_\diamond^M.$$

- Suppose that $(u, U) \in H^1(\Omega) \times \mathbb{R}^M$ satisfies the system of equations (2.2). Then for any constant $c \in \mathbb{R}$, $(u + c, U + c[1, \ldots, 1])$ also satisfies (2.2). This is natural as from a physics point of view, only potential differences can be measured. Therefore we work in the quotient space $\mathscr{H}^1$ to identify all elements of $H^1(\Omega) \times \mathbb{R}^M$ that differ up to an additive constant. As we shall see, this restores uniqueness.
- In the spirit of [28], the limit conditions $(2.2_2)$ and $(2.2_3)$ can be written at once: indeed, denoting by $\mathbb{1}_{E_m}$ the operator that maps a function $\xi \in L^\infty(E_m)$ to its extension by zero on $\partial\Omega$, and defining the operator

$$\Xi : (\xi_1, \ldots, \xi_m) \in \prod_{m=1}^{M} L^\infty(E_m) \mapsto \sum_{m=1}^{M} \mathbb{1}_{E_m} \xi_m \in L^\infty(\partial\Omega),$$

then $(2.2_2)$ and $(2.2_3)$ are equivalent to

$$\sigma \nabla u \cdot \nu + \Xi(\xi_1, \ldots, \xi_M)(u - \Xi(U_1, \ldots, U_M)) = g \text{ on } \partial\Omega.$$

- In practical applications in Electrical Impedance Tomography, the source terms $f$ and $g$ are both null, in which case the compatibility condition (2.3) becomes the usual condition $I \in \mathbb{R}_\diamond^M$. Nevertheless, non-null source terms will be used in our study in the numerical tests, to ensure smoothness on $u$ and catch the convergence rates of our numerical scheme (see subsection "Regularity results" below).
- In problem (2.2), the current configuration $I$ is imposed, and the potential configuration $U$ is measured. From a mathematical point of view, it would be equivalent to impose $U$, and measure $I$. To apply Immersed Boundary Method to our problem, it would be simpler to impose $U$, as the limit conditions are then standard. Nevertheless, in practice, it is always the current that is imposed and the potential that is measured, which breaks the symmetry as the noise is greater on the measured data then on the imposed one. Therefore we focus our study on problem (2.2).

*Well-posedness.* We now prove that problem (2.2) is well-posed. We follow the standard steps used in the case $f = 0$ and $g = 0$, based on the variational formulation of problem (2.2), that we adapt to our context.

**Proposition 2.2.** *Suppose $f$, $I$ and $g$ satisfy* (2.3). *Then problem* (2.2) *admits an unique solution $(u, U) \in \mathscr{H}^1$.*

*Proof.* Consider the variational problem: find $(u, U) \in \mathscr{H}^1$ such that, for all $(v, V)$ in $\mathscr{H}^1$,

(2.4) $$B\big((u, U), (v, V)\big) = L\big((v, V)\big),$$

with

$$B\big((u, U), (v, V)\big) = \int_\Omega \sigma \nabla u \cdot \nabla v \, dx + \sum_{m=1}^{M} \int_{E_m} \xi_m (u - U_m)(v - V_m) ds(x),$$

and

$$L\big((v, V)\big) = \int_\Omega f v \, dx + \sum_{m=1}^{M} \int_{E_m} g(v - V_m) \, ds(x) + \int_{E_c} g v \, ds(x) + I \cdot V.$$

It is clear that $B$ is bilinear and continuous on $\mathscr{H}^1 \times \mathscr{H}^1$. By definition of $\mathscr{Z}_m$, there exists $c > 0$ such that for all $m$ in $\{1, \ldots, M\}$, there exists $e_m \subset E_m$ of positive Lebesgue measure such that $\xi_m \geq c$ on $e_m$. This property combined with the assumption on $\sigma$ and Lemma 2.1 easily implies the coercivity of $B$. On the other hand, as for all $(v,V) \in H^1(\Omega) \times \mathbb{R}^M$ and all constant $c$, we have due to the compatibility condition (2.3)

$$L((v,V)) = L((v-c, V - c[1,\ldots,1])),$$

$L$ is a continuous linear form on $\mathscr{H}^1$. Therefore, Lax-Migram theorem [6, Corollary 5.8] implies the existence and uniqueness of $(u,U) \in \mathscr{H}^1$ satisfying (2.4).

Remains to prove that (2.4) is equivalent to (2.2). The fact that (2.2) implies (2.4) is standard, so we focus on the reverse implication, and prove that $(u,U) \in \mathscr{H}^1$ solution of (2.4) satisfies (2.2). Choosing $(v,V) \sim (\varphi, 0_{\mathbb{R}^M})$ in (2.4), with $\varphi \in C_c^\infty(\Omega)$, immediately implies that $-\nabla \cdot (\sigma \nabla u) = f$ in $\Omega$. Choosing then $(v,V) \sim (v_0, 0_{\mathbb{R}^M})$ with $v_0 \in H^1(\Omega)$, we obtain

$$\int_\Omega (\nabla \cdot \sigma \nabla u \, v_0 + \sigma \nabla u \cdot \nabla v_0) \, dx = \sum_{m=1}^M \int_{E_m} (g - \xi_m(u - U_m)) \, v_0 \, ds(x) + \int_{E_c} g \, v_0 \, ds(x),$$

which implies by surjectivity of the trace operator $v_0 \in H^1(\Omega) \mapsto v_{0|\partial\Omega} \in H^{1/2}(\partial\Omega)$ and the definition of the conormal derivative $\sigma \nabla u \cdot \nu$ that $u$ satisfies $(2.2)_2$ and $(2.2)_3$. In particular $\sigma \nabla u \cdot \nu$ belongs to $L^2(\partial \Omega)$.

Finally, choosing $(v,V) \sim (0, V_0)$ with $V_0 \in \mathbb{R}^M$ yields

$$I \cdot V_0 = \sum_{m=1}^M \int_{E_m} (g - \xi_m(u - U_m)) \, V_{0,m} \, ds(x) = \sum_{m=1}^M \int_{E_m} \sigma \nabla u \cdot \nu \, V_{0,m} \, ds(x) = W \cdot V_0,$$

with $W \in \mathbb{R}^M$ defined by $W_m = \int_{E_m} \sigma \nabla u \cdot \nu \, ds(x)$, which ends the proof. $\square$

*Regularity results.* The convergence results we obtain for the numerical scheme we develop in the next section require smooth solutions to our problem of interest. Unfortunately, in lots of situations, and in particular in the relevant ones for application, the solutions of (2.2) fail to be smooth.

Indeed, suppose that $\sigma$ is smooth. Due to the mixed boundary conditions $(2.2)_2$, $(2.2)_3$ and $(2.2)_4$, the potential $u$ satisfying system (2.2) might not be as smooth as in standard elliptic problem.

More precisely, for the standard Complete Electrode Model in which $\xi_m$ is a positive constant, and in practical settings where the source terms $f$ and $g$ are null, it is known that $(u,U)$ solution of (2.2) belongs to $\mathscr{H}^{2-\varepsilon}$ for all $\varepsilon > 0$, but fails to be in $\mathscr{H}^2$ for all input currents $I$ except for the null one [17, 28].

The smoothened Complete Electrode Model has been introduced in [28] precisely to overcome this difficulty. It consists in replacing $\xi_m$ constant by $\xi_m \in \mathscr{Z}_m$, compactly supported in $E_m$ and smooth. Then standard regularity results apply, leading to a smooth potential $u$. As in practice the surface admittivity is unknown, and there is no clear evidence supporting the idea that it is constant, such change of model is reasonable. Obviously, its main drawback is more complex parametrization of the contact admittivities. Our convergence results always apply to the smoothened Complete Electrode Model.

Nevertheless, regardless of the choice of $\xi_m$, there is no explicit solution for problem (2.2) when $f = g = 0$, except in very particular geometric configurations. This is why we introduce the source terms $f$ and $g$: they allow us to construct explicit solutions to problem (2.2) to numerically test our convergence results.

2.3. **An alternative formulation well-posed in $H^1(\Omega) \times \mathbb{R}^d$.** As seen previously, problem (2.2) is well posed in $\mathscr{H}^1$. In other words, the solution of (2.2) is defined up to an additive constant. From a numerical point of view, one has to fix that constant (or in other words, choose the ground level), which is usually done by adding a constraint on $U$, such that $U_1 = 0$ or $U \in \mathbb{R}_\diamond^M$: the solution is then searched in $H^1(\Omega) \times \tilde{\mathbb{R}}$, with $\tilde{\mathbb{R}} = \{U \in \mathbb{R}^m, \ U_1 = 0\}$ or $\tilde{\mathbb{R}} = \mathbb{R}_\diamond^M$, restoring uniqueness. Theoretically, the constraint can be imposed on the potential $u$ (as an example, one could impose $u$ to be mean free over $\Omega$), but this adds computational issues when solving numerically the problem.

In this study, we propose another approach, consisting in modifying problem (2.2) in a new problem, which will be well-posed in $H^1(\Omega)\times\mathbb{R}^d$, and which solution is one of the solutions to problem (2.2). Note that by construction no additional constraint is imposed on the new problem to obtain well-posedness.

*Setting $U_1 = 0$.* Let $\varepsilon > 0$ be fixed. We consider the modified problem

$$
(2.5) \quad \begin{cases}
-\nabla \cdot (\sigma \nabla u) = f & \text{in } \Omega \\
\sigma \nabla u \cdot \nu = g & \text{on } E_c \\
\sigma \nabla u \cdot \nu + \xi_m (u - U_m) = g & \text{on } E_m, \ m = 1, \ldots, M, \\
\int_{E_m} \sigma \nabla u \cdot \nu \, ds(x) + \varepsilon \delta_{m1} U_m = I_m, & m = 1, \ldots, M.
\end{cases}
$$

The boundary equation $(2.5)_4$ is equivalent to

$$\int_{E_1} \sigma \nabla u \cdot \nu \, ds(x) + \varepsilon U_1 = I_1, \quad \int_{E_m} \sigma \nabla u \cdot \nu \, ds(x) = I_m, \ \forall m \in \{2, \ldots, M\}.$$

In other words, only one equation of the initial problem (2.2) is modified.

We claim the following

**Proposition 2.3.** *Problem* (2.5) *admits an unique solution* $(u, U) \in H^1(\Omega) \times \mathbb{R}^d$. *Furthermore, if $I$, $f$ and $g$ satisfy the compatibility condition* (2.3), $U_1 = 0$ *and* $(u, U)$, *seen as an element of $\mathscr{H}^1$, is the solution of* (2.2).

*Proof.* We endow $H^1(\Omega) \times \mathbb{R}^M$ with its standard scalar product, with corresponding norm

$$\|(u, U)\|^2_{H^1(\Omega) \times \mathbb{R}^M} = \|u\|^2_{H^1(\Omega)} + \|U\|^2_{\mathbb{R}^M}.$$

First, adapting the proof of Proposition 2.2, it is readily seen that problem (2.5) is equivalent to the variationnal formulation

$$(2.6) \qquad B_1\big((u, U), (v, V)\big) = L\big((v, V)\big), \quad \forall (v, V) \in H^1(\Omega) \times \mathbb{R}^M.$$

with

$$B_1\big((u, U), (v, V)\big) = \int_\Omega \sigma \nabla u \cdot \nabla v \, dx + \sum_{m=1}^M \int_{E_m} \xi_m (u - U_m)(v - V_m) ds(x) + \varepsilon\, U_1 V_1$$

and

$$L\big((v, V)\big) = \int_\Omega f\, v \, dx + \sum_{m=1}^M \int_{E_m} g\,(v - V_m)\, ds(x) + \int_{E_c} g\, v \, ds(x) + I \cdot V.$$

Similarly, a minor adaptation of the proof of Lemma 2.1 shows that

$$\|(u, U)\|^2_{b_1} = \int_\Omega |\nabla u|^2 dx + \sum_{m=1}^M |u - U_m|^2 ds(x) + U_1^2,$$

is a norm on $H^1(\Omega) \times \mathbb{R}^M$ equivalent to the standard one. The assumptions on $\sigma$ and $\xi_m$ immediately imply the coerciveness of the continuous bilinear form $B_1$, while the linear form $L$ is clearly continuous on $H^1(\Omega) \times \mathbb{R}^M$. Therefore, Lax-Milgram theorem implies the existence and uniqueness of the solution of (2.6).

We now suppose that $I$, $f$ and $g$ satisfy (2.3). Choosing $v = 1$ and $V = [1, \ldots, 1]$ leads to

$$\varepsilon\, U_1 = \int_\Omega f \, dx + \int_{E_c} g \, ds(x) + \sum_{m=1}^M I_m = 0.$$

As a consequence, $U_1 = 0$, which immediately implies that $(u, U)$ satisfies the system of equation (2.2), which ends the proof. □

We see that solving Problem (2.5), which is a very slight variation of our initial problem (2.2), allows to retrieve the unique solution of (2.2) satisfying the additional condition $U_1 = 0$, without strongly enforcing this condition. Note that this holds for any value of the positive parameter $\varepsilon$, as long as $I$, $f$ and $g$ are compatible in the sense of (2.3).

As it allows to solve (2.2) without having to fix the ground level, from this point we will always focus on problem (2.5).

*Remark* 2.4. The focus on the first electrode is completely arbitrary here. A simple renumbering of the electrodes allows to use the same formulation to enforce the zero condition on any electrode.

*Remark* 2.5. Problem (2.5) is well-posed regardless of the compatibility of the data $I$, $f$ and $g$. If the data are not compatible, then $(u, U)$ solution of (2.5) cannot solve (2.2) as the latter has no solution.

As seen in the proof of Proposition 2.3, if $(u, U)$ is the solution of (2.2), we always have

$$(2.7) \qquad \varepsilon U_1 = \int_\Omega f \, dx + \int_{E_c} g \, ds(x) + \sum_{m=1}^{M} I_m.$$

As a consequence, if we set $\tilde{I} \in \mathbb{R}^M$ defined by

$$\tilde{I}_1 = I_1 - \left( \int_\Omega f \, dx + \int_{E_c} g \, ds(x) + \sum_{m=1}^{M} I_m \right), \quad \tilde{I}_m = I_m, \ m \neq 1,$$

we see that $(u, U)$ satisfies (2.2) with the original source terms $(I, f, g)$ replaced by $(\tilde{I}, f, g)$. This second set of data obviously satisfies the compatibility condition (2.3). In a sense, problem (2.5) correct the data so that they always satisfy the compatibility condition.

In the case of compatible data, the value of parameter $\varepsilon$ has no influence on the solution of problem (2.5), and can be freely set equal to any arbitrary value, for instance one. From a numerical point of view, this is not the case. Indeed the *numerical compatibility condition* has no reason to be exactly satisfied by the finite-difference formulation, but only up to the order of the approximation, even if the continuous data satisfies (2.3). Therefore, due to equation (2.7), the value of $\varepsilon$ influences the value of $U_1$. This causes an interplay between the two terms in the left-hand side of $(2.5)_4$ that will have an influence on the amplitude of the numerical error, but not on the convergence order itself.

*Setting* $U \in \mathbb{R}_\diamond^M$. As a final remark, we note that it is equally easy to select $(u, U)$ solution of (2.2) satisfying the additional mean free value condition on $U$, that is $U \in \mathbb{R}_\diamond^M$, without strongly enforcing this condition (that is, with a problem posed in $H^1(\Omega) \times \mathbb{R}^M$ and not in $H^1(\Omega) \times \mathbb{R}_\diamond^M$.

Indeed, it suffices to solve the following problem:

$$(2.8) \qquad \begin{cases} -\nabla \cdot (\sigma \nabla u) = f & \text{in } \Omega \\ \sigma \nabla u \cdot \nu = g & \text{on } E_c \\ \sigma \nabla u \cdot \nu + \xi_m(u - U_m) = g & \text{on } E_m, \ m = 1, \ldots, M, \\ \int_{E_m} \sigma \nabla u \cdot \nu \, ds(x) + \varepsilon \langle U \rangle = I_m, & m = 1, \ldots, M. \end{cases}$$

Here we have set $\langle U \rangle = \sum_{k=1}^{M} U_k$. The corresponding variational form is

$$B_{\langle \rangle}\big((u, U), (v, V)\big) = L\big((v, V)\big), \quad \forall (v, V) \in H^1(\Omega) \times \mathbb{R}^M.$$

with

$$B_{\langle \rangle}\big((u, U), (v, V)\big) = \int_\Omega \sigma \nabla u \cdot \nabla v \, dx + \sum_{m=1}^{M} \int_{E_m} \xi_m(u - U_m)(v - V_m) ds(x) + \varepsilon \langle U \rangle \langle V \rangle,$$

the continuous linear form $L$ being unchanged. Following the same line of reasoning used in the previous section, it is not difficult to see that (2.8) is well-posed in $H^1(\Omega) \times \mathbb{R}^M$, and that if $(I, f, g)$ satisfies (2.3), $(u, U)$ solution of (2.8) satisfies the original Complete Electrode Model equations (2.2), with furthermore $\langle U \rangle = 0$.

*Remark* 2.6. From a numerical point of view, the additional term $\langle U \rangle$ seems not practical, as it is non-local and therefore tends to fill the matrix corresponding to the discretization of (2.8). However, it appears only on equation (2.8)$_4$, and affects only the term $U \in \mathbb{R}^M$, $M$ being the number of electrodes on the boundary of $\Omega$, and therefore typically less than one hundred. As the potential $u$ is usually discretized in spaces containing thousands of degrees of freedom, in practice it does not affect significantly the matrix filling.

## 2.4. Towards the immersed boundary method.

In our numerical scheme, the domain of interest $\Omega$ is immersed into a larger domain $\Omega_e$, a square in dimension 2 and a cube in dimension 3, such that $\overline{\Omega} \subset \Omega_e$. The domain $\Omega_e$ is fully discretized with a Cartesian grid, implying that some discretization points lie in $\Omega_e \setminus \overline{\Omega}$, a part of the extended domain where there is no physical problem in our context. In that situation, a first idea is to simply set the solution to zero on the degrees of freedom lying outside of $\Omega$. But this operation is cumbersome, and not appropriate in our approach, as we need smoothness in the solution on either side of the interface $\partial\Omega$ to compute quadrature formulas on the electrodes.

We therefore choose a different approach, which consists in solving an additional partial differential equation in $\Omega_e \setminus \overline{\Omega}$, similar to the problem posed in $\Omega$. As a consequence, we introduce the following additional problem

$$
(2.9) \quad \begin{cases} -\Delta u_e = 0 & \text{in } \Omega_e \setminus \overline{\Omega} \\ u_e = 0 & \text{on } \partial\Omega_e \\ u_e = u & \text{on } \partial\Omega, \end{cases}
$$

the function $u$ being the electric potential, solution of (2.5). Clearly, such problem is well-posed in $H^1(\Omega)$, and its solution depends continuously on $u$, and therefore on the EIT source terms $(I, f, g)$. Furthermore, it will automatically inherit the smoothness of $u$ thanks to standard elliptic regularity results.

As a conclusion, our goal is now to define an immersed boundary method to reconstruct $(u, U, u_e)$, with $(u, U)$ solution of (2.5) and $u_e$ solution of (2.9).

## 3. An immersed boundary scheme for CEM

In this section we focus on the two-dimensional case, and present the discretization of problems (2.5) and (2.9) with an immersed boundary method on a Cartesian grid. The domain $\Omega$ is immersed into a larger square domain $\Omega_e$, and the boundary conditions of (2.5) are taken into account by using additional variables located on the boundary of $\Omega$. We prove that the method converges with first order accuracy. The convergence proof is based on a discrete maximum principle, used to provide estimates of the coefficients of the inverse of the discretization matrix.

### 3.1. Classification of grid points.

We consider a uniform Cartesian grid defined on the square $\Omega_e$. The grid spacing is denoted $h$. Each node of this cartesian grid is called a grid point and is denoted $M_{ij} = (x_i, y_j) = (i\,h, j\,h)$. We denote by $u_{ij}^h$ the approximation of the function $u$ at the point $(x_i, y_j)$. The set of grid points located inside the domain $\Omega_e$ is denoted $\Omega_h$.

To the purpose of the discretization of the CEM we need to define additional points on the boundary $\partial\Omega$. We define the boundary point $I_{i+1/2,j} = (x_{i+1/2}, y_j)$ as the intersection of the boundary $\partial\Omega$ and the segment $[M_{ij}M_{i+1j}]$, if it exists. Similarly, the boundary point $I_{i,j+1/2} = (x_i, y_{i,j+1/2})$ is defined as the intersection of the boundary and the segment $[M_{ij}M_{ij+1}]$. At each boundary point we create an additional unknown $u_{i+1/2,j}^h$ or $u_{i,j+1/2}^h$. The set of boundary points is denoted $\delta\Omega_h$.

We say that a grid point is irregular if one of its direct neighbors is a boundary point, see Figure 1. On the contrary, grid points that are not irregular are called regular grid points. The set of irregular grid nodes is denoted $\Omega_h^*$. The set of electrode values is denoted $E_h$.

The grid or boundary points can be denoted with letters such as $P$ or $Q$ rather than with indices such as $M_{i,j}$ if it is convenient. We also denote if more convenient $x_P$ and $y_P$ the coordinates of a point $P$.

### 3.2. Numerical scheme.

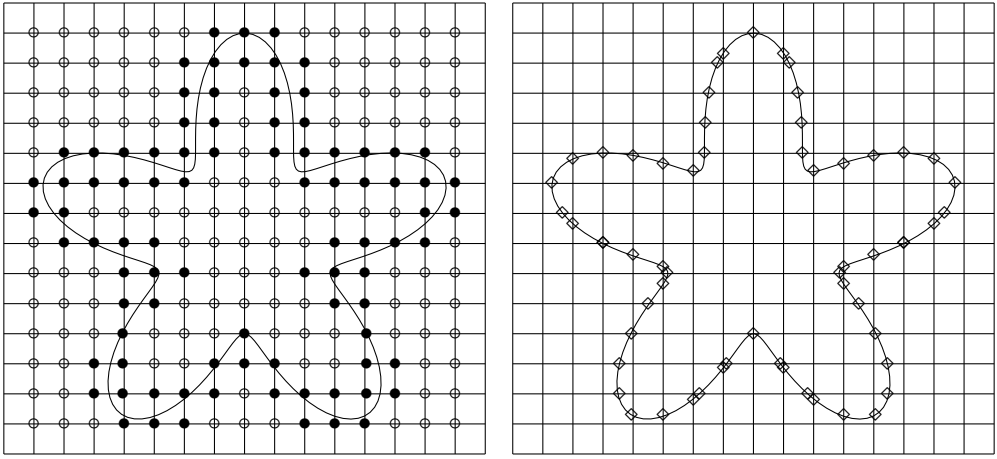

FIGURE 1. Left: regular points described by circles ∘, irregular points (belonging to $\Omega_h^*$) described by bullets •, right: boundary points (belonging to $\delta\Omega_h$).

*Elliptic operator.* To discretize the elliptic operator $(2.5)_1$ or $(2.9)_1$ on each grid point, regular or not, we use a five point stencil with the grid point $M_{i,j}$ and its nearest neighbors, boundary or grid points, in each direction. To make this explicit, we denote $u_S^h$ the value of the solution on the nearest point in the south direction, with coordinates $(x_S, y_S)$. Similarly, we define $u_N^h$, $u_W^h$ and $u_E^h$ and the associated coordinates $(x_N, y_N)$, $(x_W, y_W)$ and $(x_E, y_E)$. The discretization reads

$$
\begin{aligned}
-\Big(\nabla.(\sigma\nabla u)\Big)_{i,j}^h &= -\Big(\sigma_{i+1/2,j}\frac{u_E^h - u_{ij}^h}{x_E - x_i} - \sigma_{i-1/2,j}\frac{u_{ij}^h - u_W^h}{x_i - x_W}\Big)\frac{1}{h} \\
&\quad -\Big(\sigma_{i,j+1/2}\frac{u_N^h - u_{ij}^h}{y_N - y_j} - \sigma_{i,j-1/2}\frac{u_{ij}^h - u_S^h}{y_j - y_S}\Big)\frac{1}{h}.
\end{aligned}
$$

with $\sigma_{i+1/2,j}$ the value of $\sigma$ at point $(\frac{x_i + x_E}{2}, y_j,)$. The truncation error of this discretization is second-order accurate on regular grid points, and zeroth-order on irregular grid points.

*Flux boundary conditions.* At each boundary point we discretize the boundary conditions $(2.5)_2$ or $(2.5)_3$, depending if the boundary point belongs to an electrode or not.

For instance, on a boundary point $I_{i+1/2,j}$, the discretization reads

(3.1) $$\sigma(\nabla u \cdot \nu)_{i+1/2,j}^h + \xi_m(I_{i+1/2,j})\Big(u_{i+1/2,j}^h - U_m\Big) = 0 \text{ if } I_{i+1/2,j} \in E_m,$$

(3.2) $$\sigma(\nabla u \cdot \nu)_{i+1/2,j}^h = 0 \text{ otherwise}.$$

with $(\nabla u \cdot \nu)_{i+1/2,j}^h$ denoting the discretized normal derivative. A similar discretization is applied at each boundary point $I_{i,j+1/2}$.

The discretization of the normal derivative $(\nabla u \cdot \nu)_{i+1/2,j}^h$ depends of the local geometry of the interface. As illustrated on Figure 2, the first intersection between the normal to the boundary and the grid is located on a segment: either $[M_{i,j}, M_{i,j-1}]$, or $[M_{i,j-1}, M_{i+1,j-1}]$, or $[M_{i,j}, M_{i,j+1}]$, or $[M_{i,j+1}, M_{i+1,j+1}]$. The discrete normal derivative is computed as the normal derivative of the linear interpolant of the numerical solution on the triangle composed of the boundary point $I_{i+1/2,j}$ and the aforementioned segment.

If we denote $K$ this triangle, $(x_1, y_1)$, $(x_2, y_2)$ and $(x_3, y_3)$ its vertices, and $u_1$, $u_2$ and $u_3$ the associated values, the basis functions on the vertices for the linear interpolation write

$$\lambda_j(x,y) = \alpha_j x + \beta_j y + \gamma_j, \ j = 1, 2, 3,$$

with
$$\alpha_j = \frac{y_k - y_i}{(x_j - x_k)(y_j - y_i) - (x_j - x_i)(y_i - y_k)},$$
$$\beta_j = \frac{x_i - x_k}{(x_j - x_k)(y_j - y_i) - (x_j - x_i)(y_i - y_k)},$$
$$\gamma_j = \frac{x_k y_i - x_i y_k}{(x_j - x_k)(y_j - y_i) - (x_j - x_i)(y_i - y_k)},$$

$(\nu_x, \nu_y)$ being an approximation of the normal at the interface point. With these notations, the approximation of the normal derivative writes for instance for the interface point $I_{i+1/2,j}$

$$(\nabla u \cdot \nu)^h_{i+1/2,j} = (u_1\,\alpha_1 + u_2\,\alpha_2 + u_3\,\alpha_3)\nu_x + (u_1\,\beta_1 + u_2\,\beta_2 + u_3\,\beta_3)\nu_y.$$

This discretization is first-order accurate because it is based on a linear interpolation.

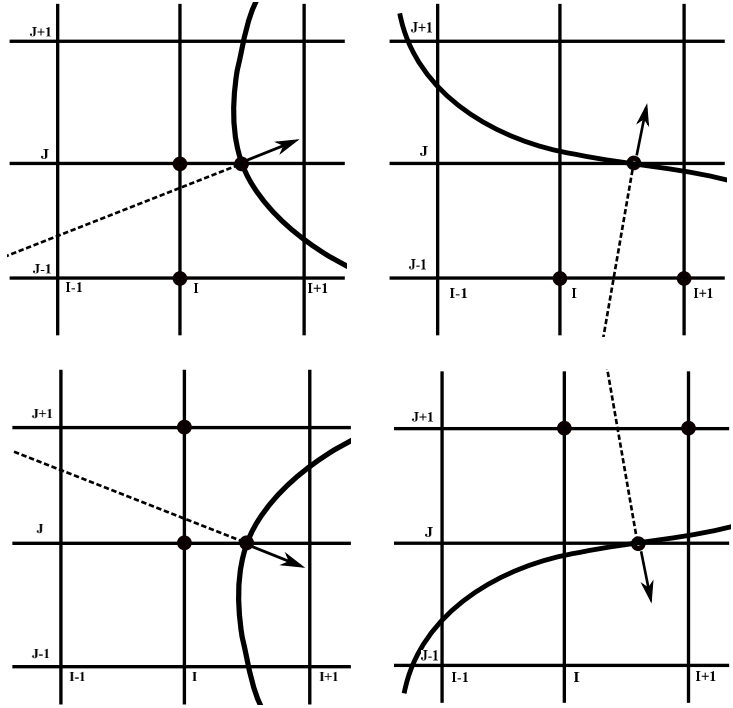

FIGURE 2. All possible stencils for the first-order flux discretization on the left side of the interface, with points involved in the discretization signaled by black circles.

*Integral boundary conditions.* We now discretize (2.5)$_4$. More precisely, using (2.5)$_3$, it is equivalent to:

$$(3.3) \qquad \int_{E_m} \bigl(g + \xi_m(U_m - u)\bigr)\,ds(x) + \varepsilon\,\delta_{m1}\,U_1 = I_m, \ m = 1,\ldots,M.$$

From now on, we focus on the equivalent boundary condition (3.3).

We discretize this integral boundary condition (3.3) on each electrode with a first-order quadrature formula based on a first-order discrete dirac function [47]. This discrete Dirac function is by construction positive and non-zero only on the irregular grid points, therefore it can be written as

$$(3.4) \qquad \sum_{P \in \Omega_h^*} \omega_P\Bigl(g(P) + \xi_m(P)(U_m - u_P)\Bigr) + \varepsilon\,\delta_{m1}\,U_1 = I_m, \ m = 1,\ldots,M.$$

with the coefficients $\omega_P$ that are the weights of the first-order quadrature formula in [47]. The truncature error of this formula is therefore first-order.

3.3. **Monotonicity of the discretization matrix.** Here we aim to prove that the discretization matrix of the linear system described before, that we denote $A_h$, is monotone. Let us first prove the following property which will be useful in the reasoning:

**Proposition 3.1.** *With the convention used for the normal to the boundary, if the minimum of $v$ is located on a boundary point, then at this boundary point the discrete normal derivative is negative. Moreover, if the approximated normal derivative at this point is zero, then the three points values involved in the stencil are equal.*

*Proof.* The approximation of the normal derivative is constant in space, because it is computed from a linear interpolation on a triangle. Therefore, if the minimum of $v$ is located on a boundary point, then at this boundary point the discrete normal derivative is negative. □

**Proposition 3.2.** *The discretization matrix of the linear system described before, that we denote $A_h$, is monotone, that is, $A_h$ is invertible and all values of $A_h^{-1}$ are non-negative.*

*Proof.* Let $v$ be an array such that all coefficients of $A_h v$ are non-negative, which we denote $A_h v \geq 0$. We aim to prove that all coefficients of $v$ are non-negative. To this purpose we consider the minimum of $v$ in the whole domain: it can be either a grid point value, or a boundary point value, or an electrode potential. We detail all these cases hereafter.

- *If the minimum is reached on a grid point:*
  In this case we denote $(i_0, j_0)$ the indices of the smallest component of $v$. We assume the grid point is a regular grid point. Otherwise the formula would have slightly different weights, but the reasoning would be the same. Using the elliptic operator inequality on this point we can write:
  $$4v_{i_0,j_0} - v_{i_0+1,j_0} - v_{i_0-1,j_0} - v_{i_0,j_0+1} - v_{i_0,j_0-1} \geq 0,$$
  and we deduce that
  $$v_{i_0+1,j_0} = v_{i_0-1,j_0} = v_{i_0,j_0+1} = v_{i_0,j_0-1} = v_{i_0,j_0}.$$

  Repeating recursively this reasoning on the neighbours of $(i_0, j_0)$, then on the neighbours of the neighbours etc, we deduce that all values of $v$ corresponding to grid or boundary points are equal to $v_{i_0,j_0}$. Therefore this case amounts to the case of the minimum located on a boundary point.

- *If the minimum is located on a boundary point not belonging to an electrode:*
  Without loss of generality, we assume that the minimum is located on $I_{i+1/2,j}$. On this interface point we have both relationships
  $$(\sigma \nabla v \cdot n)^h_{i+1/2,j} \leq 0 \text{ because the minimum is located on this point,}$$
  $$(\sigma \nabla v \cdot n)^h_{i+1/2,j} \geq 0 \text{ because all coefficients of } A_h v \text{ are non-negative,}$$
  which leads to
  $$(\sigma \nabla v \cdot n)^h_{i+1/2,j} = 0$$
  and the values of the grid points involved in the discretization of the normal derivative are also equal to the minimum. Using the reasoning of the previsous subsection, it means that all grid point and boundary point values are equal to this minimum. Therefore this case finally amounts to considering the case of the minimum on a boundary point belonging to an electrode.

- *If the minimum is located on a boundary point belonging to an electrode:*

    Without loss of generality, we assume that the minimum is located on $I_{i+1/2,j}$. On this interface point we have the relationships

    $$(\sigma \nabla v \cdot n)^h_{i+1/2,j} \leq 0,$$
    $$(\sigma \nabla v \cdot n)^h_{i+1/2,j} + \xi_m(I_{i+1/2,j})\Big(v^h_{i+1/2,j} - V_m\Big) \geq 0,$$

    the first one because the minimum is located on this point, and the second one because all coefficients of $A_h v$ are non-negative. Therefore

    $$v^h_{i+1/2,j} - V_m \geq 0,$$

    which means that

    $$V_m = v^h_{i+1/2,j},$$

    so the electrode potential $V_m$ is also the minimum of $v$. Therefore this case amounts to considering that the minimum of $v$ corresponds to an electrode potential $V_m$.

- *If the minimum corresponds to an electrode potential $V_m$ with $m \neq 1$:*

    On this electrode we have the relationship:

    $$\sum_{P \in \Omega_h^*} \omega_P \, \xi_m(P)\Big(V_m - v_P\Big) \geq 0,$$

    therefore all boundary points belonging to $E_m$ are equal to the minimum $V_m$. Using the previous reasoning it means that all values of $v$ corresponding to grid points, boundary points or electrode potentials are equal to the minimum.

- *If the minimum corresponds to the electrode potential $V_1$:*

    On this electrode we have the relationship:

    $$\sum_{P \in \Omega_h^*} \omega_P \, \xi_1(P)\Big(V_1 - v_P\Big) + \varepsilon\, V_1 \geq 0,$$

    which means that

    $$\varepsilon\, V_1 \geq 0.$$

    The minimum value of $v$ is non-negative, therefore all coefficients of $v$ are non-negative.

Therefore we have proven that if $A_h v$ is non-negative, then $v$ is also non-negative. This property has two implications: $A_h$ is invertible, and all values of $A_h^{-1}$ are non-negative. □

3.4. **Discrete Green functions.** In the following, the letters $P$ and $Q$ represent indices in the linear system, representing either discretization points (on the grid or on the boundary) or electrode values. For instance, we denote $u(P)$ the coefficient of the row of $u$ with the same index than the point $P$. Similarly, $A_h U(P)$ represents the coefficient of the $P$-th row of the array $A_h U$, and $A_h(P,Q)$ is the coefficient of the $P$-th row and $Q$-th column of the matrix $A_h$. We also define by $A_h(:,Q)$ and $A_h(P,:)$ the $Q$-th column and the $P$-th row of the matrix $A_h$.

In the spirit of [12, 51] for each $Q \in \Omega_h \cup \delta\Omega_h \cup E_h$, let us define the discrete Green's function

$$G_h(:,Q) = \Big(G_h(P,Q)\Big)_{P \in \Omega_h \cup \delta\Omega_h \cup E_h}$$

as the solution of the discrete problem:

(3.5) $$\left\{ A_h G_h(:,Q)(P) = \left\{ \begin{array}{ll} 0, & P \neq Q \\ 1, & P = Q \end{array} \right. \right.$$

The matrix $A_h$ being monotone, all values of $G_h(:,Q)$ are positive. In fact, the discrete Green functions are columns of the inverse of $A_h$. We can write the solution of the discretized problem as a sum of the source terms multiplied by the values of the discrete Green functions:

$$u_h(P) = \sum_{Q \in \Omega_h \cup \delta\Omega_h \cup E_h} G_h(P,Q) \, (A_h u_h)(Q), \quad \forall P \in \Omega_h \cup \delta\Omega_h \cup E_h.$$

We can obtain bounds on the discrete Green functions with the following result:

**Proposition 3.3.** *Let $S$ and $\tilde{S}$ be two subsets of points, $W$ a discrete function, $\alpha > 0$, $\beta > 0$ and $i$, $j \in \mathbb{N}$ such that:*

$$\begin{cases} (A_h W)(P) \geq 0, & \forall P \in \Omega_h \cup \delta\Omega_h \cup E_h \setminus \tilde{S}, \\ (A_h W)(P) \geq \alpha^{-i}, & \forall P \in S, \\ (A_h W)(P) \geq -(\beta^{-j}), & \forall P \in \tilde{S}. \end{cases}$$

*Then*

$$\sum_{Q \in S} G_h(P,Q) \leq \alpha^i W(P) + \alpha^i \beta^{-j} \sum_{Q \in \tilde{S}} G_h(P,Q), \quad \forall P \in \Omega_h \cup \delta\Omega_h \cup E_h.$$

*Proof.* Using the definition of the discrete Green functions, we can write

$$A_h W(P) \geq A_h \Big( \alpha^{-i} \sum_{Q \in S} G_h(:,Q) - \beta^{-j} \sum_{Q \in \tilde{S}} G_h(:,Q) \Big)(P), \quad \forall P \in \Omega_h \cup \delta\Omega_h \cup E_h.$$

As all coefficients of $A_h^{-1}$ are non-negative, it leads to

$$W(P) - \alpha^{-i} \sum_{Q \in S} G_h(P,Q) + \beta^{-j} \sum_{Q \in \tilde{S}} G_h(P,Q) \geq 0, \quad \forall P \in \Omega_h \cup \delta\Omega_h \cup E_h,$$

and finally we obtain the following bound:

$$\sum_{Q \in S} G_h(P,Q) \leq \alpha^i W(P) + \alpha^i \beta^{-j} \sum_{Q \in \tilde{S}} G_h(P,Q), \quad \forall P \in \Omega_h \cup \delta\Omega_h \cup E_h.$$

□

3.5. **Estimates of discrete Green functions and convergence result.** In this section, we obtain upper bounds for the discrete Green functions corresponding to the matrix arising from the discretization of (2.5)-(2.9), and deduce from them the convergence order of the numerical scheme. We make the assumption that the functions $\xi_m$ are all smooth enough so that the solution $u$ has the required regularity for the following analysis, in particular for the truncation error estimation presented in Subsection 3.2 to be valid.

**Proposition 3.4.** *For $h$ small enough, the following upper bounds hold:*

$$\sum_{Q \in \Omega_h \setminus \Omega_h^*} G_h(:,Q) \leq O(1),$$

$$\sum_{Q \in \delta\Omega_h \cup E_h} G_h(:,Q) \leq O(1),$$

$$\sum_{Q \in \Omega_h^*} G_h(:,Q) \leq O(h).$$

*Proof.* We use several discrete functions in the context of proposition 3.3 in order to obtain bounds for the different blocks of the inverse matrix $A_h^{-1}$.

- We define the array $W$ such that $W \equiv 1$. Therefore, all expressions of the linear system vanish excepted the quadrature formula involving the electrode $E_1$:

$$-\left(\nabla.(\sigma\nabla W)\right)^h_{i,j} = 0 \quad \forall M_{i,j} \in \Omega_h,$$
$$(\sigma\nabla W \cdot \nu)^h_{i+1/2,j} - \xi_m(I_{i+1/2,j})\left(W_m - W^h_{i+1/2,j}\right) = 0 \quad \forall I_{i+1/2,j} \in \delta\Omega_h,$$
$$(\sigma\nabla W \cdot \nu)^h_{i,j+1/2} - \xi_m(I_{i,j+1/2})\left(W_m - W^h_{i,j+1/2}\right) = 0, \quad \forall I_{i,j+1/2} \in \delta\Omega_h,$$
$$\sum_{P\in E_m} \omega_p \xi_m(P)\left(W_m - \tilde{W}_P\right) = 0 \quad \forall m \neq 1,$$
$$\sum_{P\in E_m} \omega_p \xi_m(P)\left(W_m - \tilde{W}_P\right) + \varepsilon W_1 = \varepsilon.$$

$G_h(:,E_1)$ representing the discrete Green function associated with the electrode $E_1$, we obtain

(3.6) $$G_h(:,E_1) = \frac{1}{\varepsilon}.$$

- We consider the exact solution $\bar{u}$ of (2.5)-(2.9), with $f=1$, $g=1$, $I_m=1$ $\forall m \neq 1$, and $I_1$ such that the compatibility condition (2.3) is satisfied. We define the array $\bar{W}$ as the discretisation of $\bar{u}$ on the grid points, boundary points and electrode values.

The discretization of the elliptic operator and the fluxes is consistent at least with first-order accuracy, excepted for the grid points in $\Omega_h^*$ where it is only zeroth-order accurate. Thus for $h$ small enough, we can write that

$$-\left(\nabla.(\sigma\bar{W})\right)^h_{i,j} \geq \frac{1}{2}, \quad \forall M_{i,j} \in \Omega_h \setminus \Omega_h^*,$$
$$-\left(\nabla.(\sigma\bar{W})\right)^h_{i,j} \geq -C_1, \quad \forall M_{i,j} \in \Omega_h^*,$$
$$(\sigma\nabla\bar{W} \cdot \nu)^h_{i+1/2,j} - \xi_m(I_{i+1/2,j})\left(\bar{W}_m - \bar{W}_{i+1/2,j}\right) \geq \frac{1}{2}, \quad \forall I_{i+1/2,j} \in \delta\Omega_h,$$
$$(\sigma\nabla\bar{W} \cdot \nu)^h_{i,j+1/2} - \xi_m(I_{i,j+1/2})\left(\bar{W}_m - \bar{W}_{i,j+1/2}\right) \geq \frac{1}{2}, \quad \forall I_{i,j+1/2} \in \delta\Omega_h$$
$$\sum_{P\in E_m} \omega_p \xi_m(P)\left(\bar{W}_m - \bar{W}_P\right) \geq \frac{1}{2}, \quad \forall m \neq 1$$
$$\sum_{P\in E_1} \omega_p \xi_m(P)\left(\bar{W}_1 - \bar{W}_P\right) + \varepsilon\bar{W}_1 \geq -2|I_1|,$$

with $C_1$ a strictly positive constant. It can also be re-written as:

$$(A_h\bar{W})(P) \geq \frac{1}{2}, \quad \forall P \in (\Omega_h \setminus \Omega_h^*) \cup (E_h \setminus E_1),$$
$$(A_h\bar{W})(P) \geq \frac{1}{2}, \quad \forall P \in \delta\Omega_h,$$
$$(A_h\bar{W})(P) \geq -C_1, \quad \forall P \in \Omega_h^*,$$
$$(A_h\bar{W})(E_1) \geq -2|I_1|,$$

and it leads to:

(3.7) $$\sum_{Q\in(\Omega_h\setminus\Omega_h^*)\cup(E_h\setminus E_1)} G_h(:,Q) + \sum_{Q\in\delta\Omega_h} G_h(:,Q) \leq 2\bar{W} + 2C_1 \sum_{Q\in\Omega_h^*} G_h(:,Q) + 4|I_1|G_h(:,E_1).$$

- We define the array $\tilde{W}$ such that $\tilde{W} = 1$ for all points in $\Omega_h$, and $\tilde{W} = 0$ for all points in $\delta\Omega_h \cup E_h$. Because of the jump in the values of $\tilde{W}$ between $\Omega_h$ and $\delta\Omega_h \cup E_h$, when we apply the discretization to $\tilde{W}$ we obtain

$$-\left(\nabla.(\sigma\nabla\tilde{W})\right)^h_{i,j} = 0 \quad \forall M_{i,j} \in \Omega_h \setminus \Omega_h^*,$$

$$-\left(\nabla.(\sigma\nabla\tilde{W})\right)^h_{i,j} \geq \frac{C_2}{h^2} \quad \forall M_{i,j} \in \Omega_h^*,$$

$$(\sigma\nabla\tilde{W}\cdot\nu)^h_{i+1/2,j} - \xi_m(I_{i+1/2,j})\left(\tilde{W}_m - \tilde{W}_{i+1/2,j}\right) \geq -\frac{C_3}{h} \quad \forall I_{i+1/2,j} \in \delta\Omega_h,$$

$$(\sigma\nabla\tilde{W}\cdot\nu)^h_{i,j+1/2} - \xi_m(I_{i,j+1/2})\left(\tilde{W}_m - \tilde{W}_{i,j+1/2}\right) \geq -\frac{C_3}{h}, \quad \forall I_{i,j+1/2} \in \delta\Omega_h$$

$$\sum_{P\in E_m} \omega_p\, \xi_m(P)\left(\tilde{W}_m - \tilde{W}_P\right) + \varepsilon\delta_{m1}\tilde{W}_1 = 0 \quad \forall m,$$

with $C_2$ and $C_3$ trictly positive constants. We thus obtain

$$\frac{C_2}{h^2}\sum_{Q\in\Omega_h^*} G_h(:,Q) \leq \tilde{W} + \frac{C_3}{h}\sum_{Q\in\delta\Omega_h} G_h(:,Q),$$

and then:

$$\sum_{Q\in\Omega_h^*} G_h(:,Q) \leq \frac{h^2}{C_2}\tilde{W} + h\frac{C_3}{C_2}\sum_{Q\in\delta\Omega_h} G_h(:,Q). \tag{3.8}$$

Combining 3.8 with 3.7 we obtain

$$\sum_{Q\in(\Omega_h\setminus\Omega_h^*)\cup(E_h\setminus E_1)} G_h(:,Q) + \sum_{Q\in\delta\Omega_h} G_h(:,Q) \leq 2\bar{W} + 2C_1\frac{h^2}{C_2}\tilde{W} + 4|I_1|G_h(:,E_1) \tag{3.9}$$

$$+ 2C_1\frac{C_3}{C_2}h\sum_{Q\in\delta\Omega_h} G_h(:,Q), \tag{3.10}$$

and for $h$ small enough it leads to

$$\sum_{Q\in(\Omega_h\setminus\Omega_h^*)\cup(E_h\setminus E_1)} G_h(:,Q) + \sum_{Q\in\delta\Omega_h} G_h(:,Q) \leq \bar{C} \tag{3.11}$$

with $\bar{C}$ a positive constant. Re-injecting this inequality into 3.8 we obtain

$$\sum_{Q\in\Omega_h^*} G_h(:,Q) \leq \frac{h^2}{C_2}\tilde{W} + h\frac{C_3}{C_2}\bar{C}. \tag{3.12}$$

The inequalities 3.6, 3.11 and 3.12 provide the bounds of the proposition.

$\square$

**Proposition 3.5.** *If we denote $\bar{u}$ the exact solution of (2.5)-(2.9) and $u_h$ the numerical solution, the local error $|\bar{u}(P) - u_h(P)|$ satisfies for $h$ small enough:*

$$|\bar{u}(P) - u_h(P)| \leq O(h). \tag{3.13}$$

*Therefore the numerical solution converges with first-order accuracy to the exact solution in $L^\infty$-norm.*

*Proof.* We denote $\tau$ the truncation error array arising from the discretization of (2.5)-(2.9). We have

$$|\bar{u}(P) - u_h(P)| = \left|\sum_{Q\in\Omega_h\cup\delta\Omega_h\cup E_h} G_h(P,Q)\tau(Q)\right|$$

$$\leq \left|\sum_{Q\in\Omega_h^*} G_h(P,Q)\tau(Q)\right| + \left|\sum_{Q\in\Omega_h\setminus\Omega_h^*} G_h(P,Q)\tau(Q)\right| + \left|\sum_{Q\in\delta\Omega_h\cup E_h} G_h(P,Q)\tau(Q)\right|,$$

which, according to the previous results, immediately implies

$$|\bar{u}(P) - u_h(P)| \leq O(1)\left|\sum_{Q\in\Omega_h^*} G_h(P,Q)\right| + O(h^2)\left|\sum_{Q\in\Omega_h\setminus\Omega_h^*} G_h(P,Q)\right| + O(h)\left|\sum_{Q\in\delta\Omega_h\cup E_h} G_h(P,Q)\right|$$

$$\leq O(1)O(h) + O(h^2)O(1) + O(h)O(1) = O(h).$$

□

## 4. Numerical experiments

We now study the efficiency of the immersed boundary method for Electrical Impedance Tomography introduced in Section 3.2. More precisely, we test the method on three different two-dimensional problems:

- the direct problem of EIT, which in our context reads: for given source terms $f$, $g$ and $I$, compute $(u, U)$ solution of (2.5) (in the sense of Proposition 2.3). In that case, the geometry (i.e. $\Omega$ and the electrodes' position) and the conductivity are supposed to be known parameters. This is a well-posed problem, and our numerical method is designed to obtain a good reconstruction of its solution without accurate discretization of the geometry.
- the inverse problem of EIT: for given input currents $I \in \mathbb{R}_\diamond^M$, and corresponding measured voltages $U$, reconstruct the interior conductivity $\sigma$. In other word, knowing the couple $(I, U)$, we search for $\sigma$ such that $(u, U) = (u(\sigma), U(\sigma))$ solves (2.5) with source terms $f = 0$, $g = 0$ and $I$. Here again, the geometry is supposed to be known. This of course is the central problem in Electrical Impedance Tomography, and therefore is the main validation of our method.
- a geometric inverse problem in EIT: for given input currents $I \in \mathbb{R}_\diamond^M$, and corresponding measured voltages $U$, reconstruct the positions of the electrodes *knowing $\Omega$ and the conductivity $\sigma$*. Obviously, this problem does not make much sense from a practical point of view. But on one hand it can be seen at a first step to tackle the real problem of imperfectly known geometrical setting in EIT, which is known to severely deteriorate the quality of most conductivity's reconstruction algorithms (at the exception of the ones that are specifically designed to cope with this problem, see [2, 15, 33, 34] and the references therein). On the second hand, it is a problem of interest in our context, as the electrodes are not exactly approximated in the immersed boundary method we proposed. Therefore, knowing the high-sensibility of EIT-measurements to electrodes' position, it is a challenging problem.

4.1. **Representation of the geometry.** In the following experiments, the domain $\Omega$ is represented by a polar parametrization of its boundary $\partial\Omega$. More precisely, for a given vector $\alpha = [\alpha_0, \alpha_1, \ldots, \alpha_{2N}] \in \mathbb{R}^{2N+1}$, we set

$$\partial\Omega = \{r(\theta)\mathbf{u}(\theta), \theta \in [0, 2\pi]\},$$

with $\mathbf{u}(\theta) = [\cos(\theta), \sin(\theta)]$, and

$$r(\theta) = \alpha_0 + \sum_{k=1}^{N} (\alpha_k \cos(k\theta) + \alpha_{k+N}\sin(k\theta)).$$

The shape parameters $[\alpha_0, \ldots, \alpha_{2N}]$ are chosen so that $0 < r(\theta) < 2$ for all $\theta \in [0, 2\pi]$. Finally, we set

$$\Omega = \left\{\rho[\cos(\theta), \sin(\theta)] \in \mathbb{R}^2, \theta \in [0, 2\pi],\ 0 \leq \rho < r(\theta_0)\right\}.$$

Note that, with such parameterization of $\Omega$, the outward unit vector at a point $x = \|x\|[\cos(\theta_x), \sin(\theta_x)] \in \partial\Omega$, denoted $\nu(x)$, is given by

$$\nu(x) = \frac{1}{\rho(\theta_x)} \left[r(\theta_x)\mathbf{u}(\theta_x) - r'(\theta_x)\mathbf{v}(\theta_x)\right],$$

with $\mathbf{v}(\theta) = [-\sin(\theta), \cos(\theta)]$ and $\rho(\theta) = \sqrt{r(\theta)^2 + r'(\theta)^2}$. We also set

$$\tau(x) = \frac{1}{\rho(\theta_x)} \left[r'(\theta_x)\mathbf{u}(\theta_x) + r(\theta_x)\mathbf{v}(\theta_x)\right],$$

the vector tangential to $\partial\Omega$ at the point $x$, such that $(\nu(x), \tau(x))$ is a direct basis of $\mathbb{R}^2$.

In the following tests, we use three different geometries for $\Omega$: $\Omega_1$ is a disk of radius 1.5, corresponding to $\alpha = [1.5]$, $\Omega_2$ corresponds to the choice $\alpha = [1.51, 0.01, 0.05, 0.2, 0.035, 0.01, 0.1]$, and $\Omega_3$ corresponds to the choice $\alpha = [1.6, 0.002, 0.01, 0.003, 0.035, 0.2, 0.15]$.

The $M$-electrodes $E_1, \ldots, E_M$ are parameterized by two vectors $(\Theta^1, \Theta^2) \in \mathbb{R}^M \times \mathbb{R}^M$ verifying

$$\Theta^1_1 < \Theta^2_1 < \Theta^1_2 < \ldots < \Theta^2_M < \Theta^1_1 + 2\pi,$$

such that

$$E_m = \left\{r(\theta)\mathbf{u}(\theta),\ \theta \in [\Theta^1_m, \Theta^2_m]\right\}.$$

As in practical application the length of the electrodes is prescribed, we usually set $\Theta^1$, choose a length size $L > 0$, and compute $\Theta_2$ by solving the equation

$$(4.1) \qquad L = \int_{\Theta^1_k}^{\Theta^2_k} \rho(\theta)\,d\theta.$$

Unless otherwise specified, in the experiments below we use 16 electrodes, defined by $\Theta^1_k = -\pi + (k-1)\frac{\pi}{8}$ and $L = 0.35$.

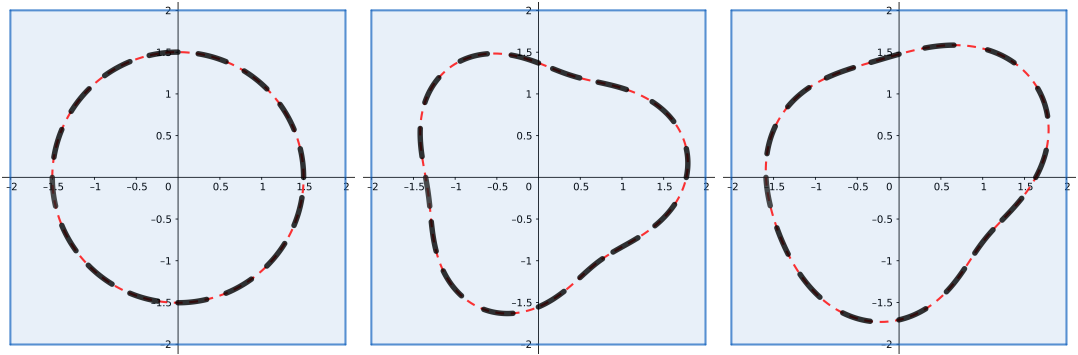

FIGURE 3. The domains $\Omega_1$, $\Omega_2$ and $\Omega_3$ (in red) in the reference square $[-2, 2] \times [-2, 2]$ (in blue), with our usual configuration of 16 electrodes (in bold black)

*Remark* 4.1. In this article we use a specific polar parametrization of the interface, which allows us to compute explicitly all normals to the interface. However, such parametrization implies that the domain of interest is star-shaped, which might not be the case in certain applications.

In such situations, it is usually still possible to use a description of the interface with a level-set function. Such a description amounts to represent implicitly the interface as the zero level set of a function. The level set method was introduced by Osher and Sethian [44]. We refer the interested reader to [43, 45, 46].Although any smooth enough function cancelling at the interface location can be used, in practice the *signed distance* to the interface is very often used:

$$\varphi(x) = \begin{cases} \text{dist}_{\partial\Omega}(x) & \text{outside of the interface } \partial\Omega, \\ 0 & \text{on the interface,} \\ -\text{dist}_{\partial\Omega}(x) & \text{inside of the interface } \partial\Omega. \end{cases} \tag{4.2}$$

We recall here some properties occurring in this case:
- As recalled in [14], the level-set, being a distance function, is 1-Lipshitz and almost everywhere differentiable. Moreover, if $\phi$ is differentiable at a point $x$, then it satisfies the so-called Eikonal equation at $x$:
$$\|\nabla\varphi(x)\| = 1.$$
- The smoothness of the level-set is in fact strongly related to the smoothness of the interface: as proved in [25, p.355], if the interface is $C^2$, then there exists a real $r_0 > 0$ such that the level-set is $C^2$ for all $x, y$ such that $|\varphi(x, y)| < r_0$.
- The outward normal vector of the isoline of $\varphi$ passing on $x$, denoted $\nu(x)$, can be expressed where $\varphi$ is differentiable as
$$\nu(x) = \frac{\nabla\varphi(x)}{\|\nabla\varphi(x)\|}, \tag{4.3}$$

4.2. **Direct problems.** First of all, we investigate the ability of the proposed immersed boundary method to efficiently compute solutions of the Electrical Impedance Tomography problem (2.5). We highlight that the convergence results announced in Section 3.5 are indeed observed in practice for solutions smooth enough, and present several examples of reconstruction of the inner potential.

4.2.1. *Convergence of the method.* To verify that the convergence rates proved in Section 3.5 are indeed observed in practice, we proceed as follows: we set $\Omega_e = [-2, 2] \times [-2, 2]$, $u(x, y) = \sin(x\, y)$, $U = 0_{\mathbb{R}^{16}}$, and compute the corresponding source terms $f$, $g$ and $I$ in (2.2). In particular, to obtain the vector of input currents $I$, we compute the corresponding integrals appearing in $(2.2)_4$ using a numerical integration scheme.

We then use the immersed boundary method described in Section 3.2 to obtain an approximate solution $u_h$ to problem (2.2) for different values of the mesh parameter $h$. To test also the convergence of the derivatives, we compute an approximation of the gradient $du_h$ by numerical differentiation.

In the following, we compute $u_h$ and $du_h$ for various values of $h$ ranging from $\frac{1}{25}$ to $\frac{1}{225}$. The parameter $\varepsilon$ in $(2.5)_4$ is chosen small, of order $\sim 10^{-10}$. The results are displayed in Figure 4: we observe a linear convergence both for $u$ and its gradient.

*Remark* 4.2. The linear convergence of $u_h$ to $u$ was expected. The linear convergence of the gradient is not surprising, as the construction of the numerical method is based on a first-order discretization of fluxes near the boundary of the domain. Nonetheless, to prove it is not obvious and non-standard within the finite difference framework. We left the study of the convergence of the gradient for future works.

*Remark* 4.3. We have presented here numerical results for a given value of the parameter $\varepsilon$. Let us remark that changing the value of the parameter $\epsilon$ does not change the convergence behaviour. However, the value of $\varepsilon$ affects the amplitude of the error *at fixed h*. This is expected since the value of $U_1$ depends on the numerical compatibility condition, and therefore is not exactly zero. Therefore, there is an interplay between the therm $\varepsilon U_1$ and the integral term in the formula $(2.5)_4$ that influences the value of both terms.

On Figure 5, we present a comparison of numerical errors for the generalized CEM, with smooth spatially varying admittivities, and for the classical CEM with $\xi_m = \frac{1}{z_m}$ constant. The geometry considered is a circle of radius $R = 0.5$ embedded in a square $[-1, 1] \times [-1, 1]$, with 4 electrodes equally spaced.The discretization parameter $h$ ranges from $\frac{1}{50}$ to $\frac{1}{250}$. The convergence rate displayed on the Figure is computed with the classical CEM. We observe that the amplitude of the numerical errors for both models as well as their rate of convergence is very close.

Finally, we present some reconstructions of the inner potential, with $\Omega_e = [-2,2] \times [-2,2]$ discretized with a $700 \times 700$ Cartesian grid, and in different geometries, with different background conductivities and various patterns applied to the electrodes, see Figures 6 and 7.

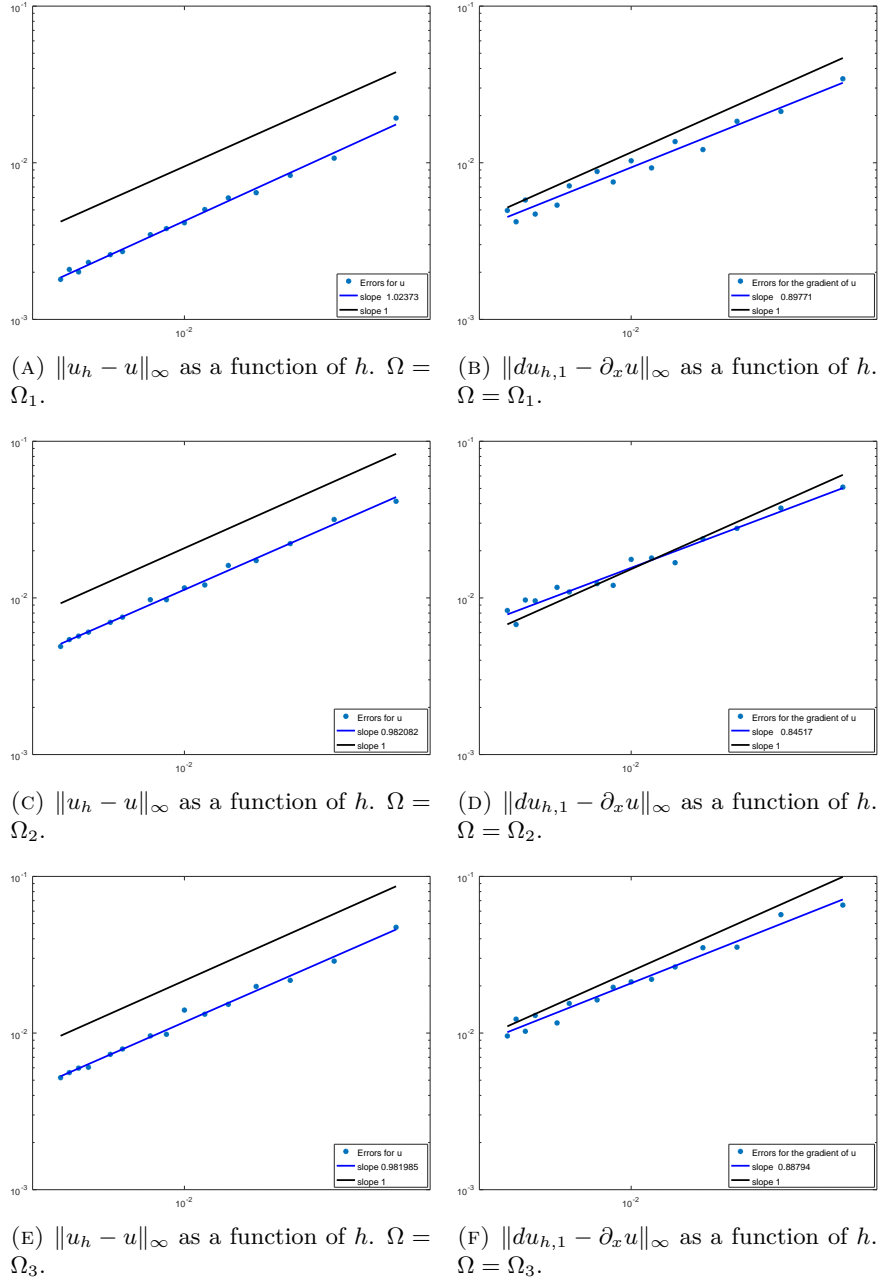

(A) $\|u_h - u\|_\infty$ as a function of $h$. $\Omega = \Omega_1$.

(B) $\|du_{h,1} - \partial_x u\|_\infty$ as a function of $h$. $\Omega = \Omega_1$.

(C) $\|u_h - u\|_\infty$ as a function of $h$. $\Omega = \Omega_2$.

(D) $\|du_{h,1} - \partial_x u\|_\infty$ as a function of $h$. $\Omega = \Omega_2$.

(E) $\|u_h - u\|_\infty$ as a function of $h$. $\Omega = \Omega_3$.

(F) $\|du_{h,1} - \partial_x u\|_\infty$ as a function of $h$. $\Omega = \Omega_3$.

FIGURE 4. Convergence study with manufactured solution $u(x,y) = \sin(x\,y)$ and $U = 0_{\mathbb{R}^{16}}$ in each domain.

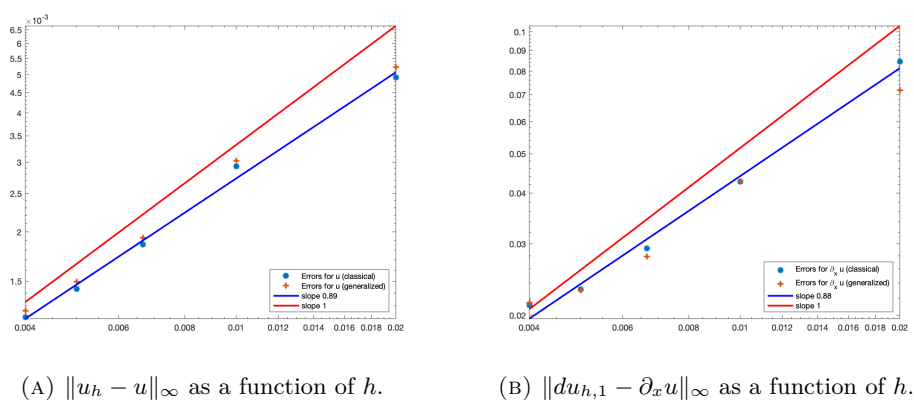

(A) $\|u_h - u\|_\infty$ as a function of $h$.

(B) $\|du_{h,1} - \partial_x u\|_\infty$ as a function of $h$.

FIGURE 5. Convergence study with manufactured solution $u(x,y) = \exp(x^2+y^2)$ and $U = 0.5\exp(R^2)1_{\mathbb{R}^4}$. Blue points stand for the numerical errors of classical CEM, red crosses for the numerical errors of generalized CEM

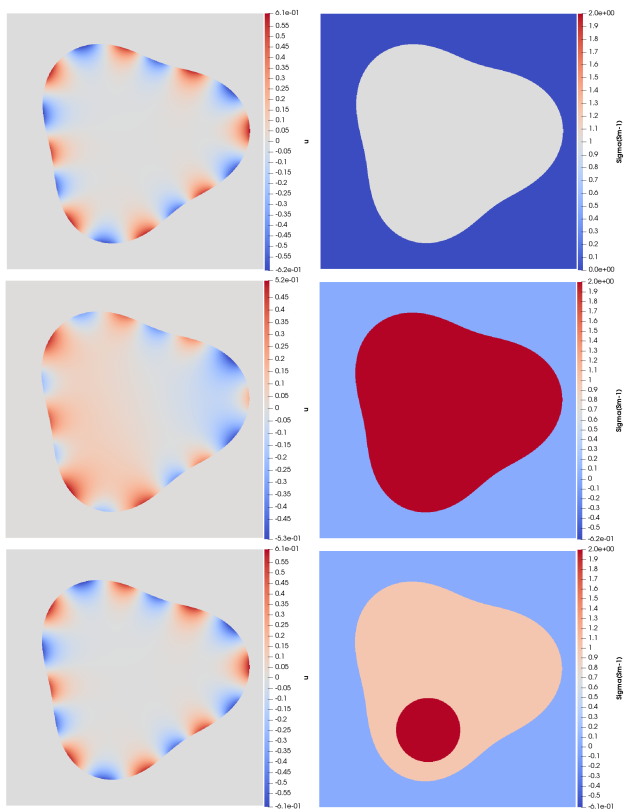

FIGURE 6. Resolution of Problem (2.5) in domain $\Omega_2$ with input current $I_m = (-1)^m$. Left: inner potential. Right: conductivity.

4.3. **Inverse problems.** From now on, we consider the standard Complete Electrode Model. In other words, in Problems (2.2) and (2.5), we set $\xi_m = \frac{1}{z_m}$, $z_m$ being the constant and positive contact impedance on each electrodes.

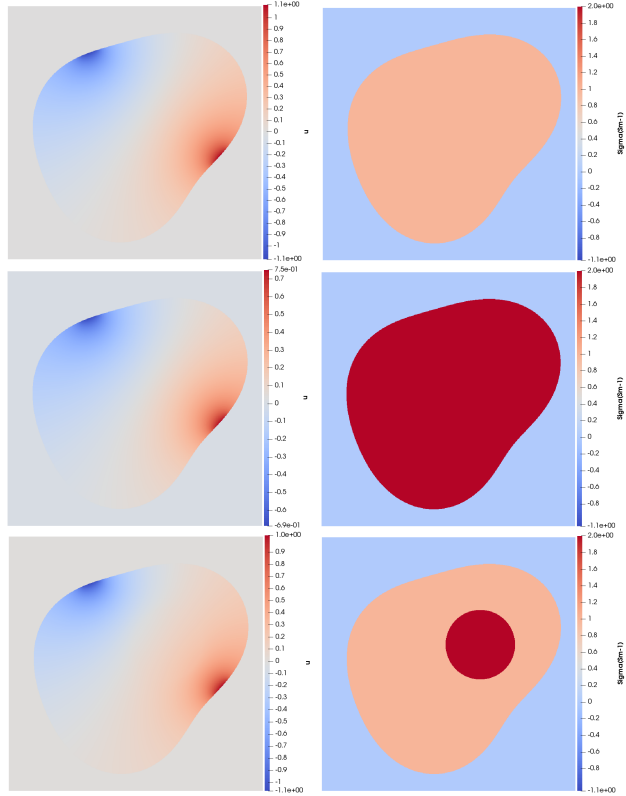

FIGURE 7. Resolution of Problem (2.5) in domain $\Omega_3$ with input current $I_m = \delta_{m\,8} - \delta_{m\,14}$. Left: inner potential. Right: conductivity.

4.3.1. *Reconstruction of the conductivity, knowing the geometry.* We now focus on the reconstruction of the conductivity $\sigma_*$ in a practical Electrical Impedance Tomography context: for a family of input currents $(I^1, \ldots, I^P) \in \mathbb{R}^M_\diamond$, $P \in \mathbb{N}_*$, the corresponding (potentially noisy) electrodes' potentials $U^1, \ldots, U^P$ are measured. They correspond to the second part of the solution of (2.2) with $\sigma = \sigma_*$ and $f = g = 0$. To fix the constant, we set $U^p_1 = 0$ for all $p$ in $\{1, \ldots, P\}$. As a consequence, the potentials are also the second part of the solution of Problem (2.5)

To simplify notations, we gather all input currents and electrodes' potentials corresponding to a certain conductivity $\sigma$ in $M \times P$ matrices

$$\mathscr{I} = \begin{bmatrix} I^1 & \ldots & I^P \end{bmatrix}, \quad \mathscr{U} = \begin{bmatrix} U^1, \ldots, U^P \end{bmatrix}.$$

The measurement matrix $\mathscr{U}$ depends linearly on the input currents $\mathscr{I}$ and non-linearly on the unknown conductivity $\sigma$. We will sometimes make this dependency explicit by writing $\mathscr{U}(\sigma, \mathscr{I})$. We also note

$$\mathscr{U}_{meas} \approx \mathscr{U}(\sigma_*, \mathscr{I}),$$

the matrix which would idealistically corresponds to our actual measurements, but in practice is only close to it, as in applications measurements are always corrupted by noise. The inverse conductivity problem can be formulated in the following concise way: knowing $(\mathscr{I}, \mathscr{U}_{meas})$ reconstruct $\sigma_*$.

Over time, sophisticated methods have been introduced to reconstruct either the conductivity or the support of inclusions from (continuous or electrodes' like) boundary measurements [7, 10, 20, 21, 26, 27, 29, 31, 35, 41, 50]. In our present study, we will use a comparatively very naive approach to reconstruct $\sigma_*$, that is a basic least-square minimization approach combined with Tychonov regularization. We choose this on purpose, as our primary objective is to test the ability of our *direct method* to give

IMMERSED BOUNDARY METHOD FOR CEM IN EIT 23sufficiently good approximation of the solution to problem (2.5) to be efficiently used in an *standard algorithm for the inverse conductivity problem*.

More precisely, we seek to minimize the discrepancy between the solution of the direct problem (2.5), corresponding to a guessed conductivity $\sigma$, and the measured data, i.e., we seek $\sigma$ such that the functional

$$(4.4) \qquad F(\sigma) = \|\mathscr{U}(\sigma, I) - \mathscr{U}_{meas}\|_2^2$$

is minimal. The problem being ill-posed, a regularization is needed. As a consequence, we actually minimize the functional

$$(4.5) \qquad F(\sigma) = \|\mathscr{U}(\sigma, I) - \mathscr{U}_{meas}\|_{\mathbb{R}^M}^2 + \frac{\epsilon}{2} R(\sigma)$$

where $R$ is a properly chosen regularisation functional and $\epsilon > 0$ a regularisation parameter. In other words, we use a deterministic Tychonov regularization of our problem.

From now on, we assume that $\sigma_* \in L^\infty \cap H^1(\Omega)$. We also assume that there exists *a known conductivity* $\sigma_\star \in L^\infty \cap H^1(\Omega)$ such that $\sigma_* - \sigma_\star \in H_0^1(\Omega)$. This is typically the case in applications where one searches perturbations of a *known* background conductivity $\sigma_\star$, assuming that the perturbations do not touch the boundary of the domain (in that situation, the pertubations are $\sigma_* - \sigma_\star$). It is then natural to seek to minimize the following functional

$$(4.6) \qquad F(\sigma) = \frac{1}{2}\|\mathscr{U}(\sigma, I) - \mathscr{U}_{meas}\|_{\mathbb{R}^M}^2 + \frac{\epsilon}{2}\|\sigma - \sigma_\star\|_{H^1(\Omega)}^2.$$

The strategy of minimization is simple: for a conductivity $\sigma$ and a current pattern $I \in \mathbb{R}_\diamond^M$, we denote $(u(\sigma, I), U(\sigma, I))$ the solution of (2.5) with $f = g = 0$. Then, for a given conductivity $\sigma_n$, define $(d_x u_n, d_y u_n)$ by

$$d_x u_n = \left[\frac{\partial u(\sigma_n, I^1)}{\partial x}, \ldots, \frac{\partial u(\sigma_n, I^P)}{\partial x}\right] \in L^2(\Omega)^d, \quad d_y u_n = \left[\frac{\partial u(\sigma_n, I^1)}{\partial y}, \ldots, \frac{\partial u(\sigma_n, I^P)}{\partial y}\right] \in L^2(\Omega)^d.$$

Now, with $\mathscr{U}_{meas} = [U_{meas}^1, \ldots, U_{meas}^P]$, we define $(d_x w_n, d_y w_n)$ by

$$d_x w_n = \left[\frac{\partial u(\sigma_n, U(\sigma_n, I^1) - U_{meas}^1)}{\partial x}, \ldots, \frac{\partial u(\sigma_n, U(\sigma_n, I^P) - U_{meas}^P)}{\partial x}\right]$$

and

$$d_y w_n = \left[\frac{\partial u(\sigma_n, U(\sigma_n, I^1) - U_{meas}^1)}{\partial y}, \ldots, \frac{\partial u(\sigma_n, U(\sigma_n, I^P) - U_{meas}^P)}{\partial y}\right].$$

Finally, we define $\delta\sigma_n$ as the unique $v \in H_0^1(\Omega)$ such that

$$-\Delta v + v = \epsilon\bigl(\Delta(\sigma_n - \sigma_\star) - (\sigma_n - \sigma_\star)\bigr) + d_x u_n^T d_x w_n + d_y u_n^T d_y w_n.$$

Then it is not difficult to prove the following result.

**Lemma 4.4.** *Suppose that $\delta\sigma_n$ belongs to $L^\infty(\Omega)$. There exists $T_n > 0$ such that for all $t \in (0, T_n)$,*

$$F(\sigma_n + t\delta\sigma_n) \leq F(\sigma_n).$$

The proof of this point is standard, and is given in the Appendix.

Our algorithm of minimization is then the following:
- Choose $\sigma_0 \in L^\infty \cap H^1(\Omega)$ (for example, $\sigma_0 = \sigma_\star$),
- Until the following stopping criteria is reached:

$$|F(\sigma_n + t_n\delta\sigma_n) - F(\sigma_n)| < \tau|\delta\sigma_n|, \text{ with } \tau \text{ small enough}$$

  (1) For a given conductivity $\sigma_n$, compute $\delta\sigma_n$ as described above,
  (2) Define $\sigma_{n+1} = \sigma_n + t_n\delta\sigma_n$ with $t_n = \arg\min_{t \in (0, T_n)} F(\sigma_n + t\delta\sigma_n)$.
- End.

Note that, at each step $n$ of this algorithm, we have to solve $P$ problems (2.5) to compute $d_x u_n$ and $d_y u_n$, then again $P$ problems (2.5) to compute $d_w w_n$ and $d_y w_n$, and a single Laplace-type problem to compute $\delta \sigma_n$, that is $2P+1$ elliptic problems. Luckily, the first $2P$ problems are exactly (2.5) with a fixed given conductivity $\sigma_n$, and only the right-hand side changing, meaning that the corresponding system to inverse is the same for the $2P$ problems. Hence the computation of $\delta \sigma_n$ is actually very fast.

Our numerical experiments show that, for $t$ close to 0, the function $t \mapsto F(\sigma_n + t\delta\sigma_n)$ behave like a quadratic function in $t$ with a single local minimum. We therefore search for $t_n$ using a *Golden-section search*.

We now present some conductivity's reconstructions in the domains presented in Figure 3. In each test, the background conductivity $\sigma_\star$ is equal to one, and

- first case: we place a circular inclusion in the center.
- second case : we place a similar inclusion closer to the boundary.
- third case : we place two inclusions of different amplitudes inside the domain.

The measurement matrix $\mathscr{U}_{meas}$ is obtained by solving (2.5) using our numerical method, in a very refined mesh and for $P = 15$ linearly independent current patterns. The mesh used to compute $\mathscr{U}_{meas}$ differs from the mesh used to reconstruct the conductivity in order to avoir inverse crime. To obtain the noisy data, we then perturb the measurements with an additive Gaussian noise, scaled so that the relative amplitude of the noise is set to $\delta \in [0,1]$. Finally, we set the regularization parameter to $\epsilon = 10^{-4}$.

Figures 8 and 9 present reconstructions in the different geometrical domains with various levels of added Gaussian noise, discretized by a $400 \times 400$ Cartesian grid. As can be seen, the reconstructions are correct even with this very naive approach, meaning that our immersed boundary methods can be efficiently used in the inverse conductivity problem.

*Remark* 4.5. In practical applications, the conductance $\xi_m$ (or equivalently the contact impedances $z_m$) are unknown, and should be also reconstructed. This can be done without additional difficulties, using that $(u, U)$ solution of (2.5) is Frechet differentiable with respect to $\xi_m$, with an explicit expression of the Frechet derivative.

4.3.2. *Reconstruction of the electrodes' positions, knowing the conductivity.* We now focus on the geometric inverse problem described in the introduction, which can be summarized as follows: *from the knowledge of $(\mathscr{I}, \mathscr{U}_{meas})$, find the locations of the electrodes, in other words the parameters $\Theta^1$ and $\Theta^2$, all the other parameters ($\alpha$, $\sigma$, ...) being known.*

As already said, even if this problem makes little sense from a practical point of view, it is interesting and challenging for our method as (2.2) is known to be very unstable with respect to the geometrical setting, and our Boundary Immersed Method is, by design, not precise in terms of geometry of $\partial \Omega$.

Adapting the notations of previous Section 4.3.1, we denote by $U = U(\Theta^1, \Theta^2, I)$ the (second term of the) solution of (2.5) with $f = g = 0$, $I \in \mathbb{R}^M_\diamond$, and electrodes' positions determined by the angles $\Theta^1$ and $\Theta^2$ (see Section 4.1 for the parametrization of the geometry). Now, for a set of input currents $(I^1, \ldots, I^P) \in (\mathbb{R}^M_\diamond)^P$, we denote

$$\mathscr{I} = \left[I^1, \ldots, I^P\right], \quad \mathscr{U} = \left[U^1, \ldots, U^P\right].$$

We also set $\mathscr{U}_{meas}$ to be the following $M \times P$ possibly noisy measurements matrix, approximating the exact measurement matrix corresponding to the exact electrodes positions $\Theta^1_*$ and $\Theta^2_*$:

$$\mathscr{U}_{meas} \approx \mathscr{U}(\Theta^1_*, \Theta^2_*, \mathscr{I}).$$

Our strategy is similar as before: we seek to minimize the functional

$$F(\Theta^1, \Theta^2) = \|\mathscr{U}(\Theta^1, \Theta^2, I) - \mathscr{U}_{meas}\|_2^2 + \frac{\epsilon}{2}\|\Theta^1 - \Theta^1_\star\|_{\mathbb{R}^M}^2 + \frac{\epsilon}{2}\|\Theta^2 - \Theta^2_\star\|_{\mathbb{R}^M}^2.$$

As previously, $\Theta^k_\star$ has to be think as a good guess for the actual location of the $k$-th electrode, around which we search the actual position of the electrode.

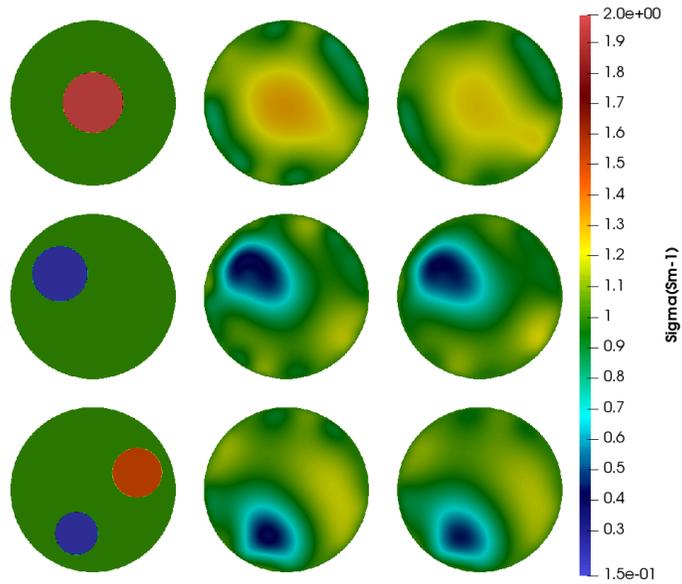

FIGURE 8. Reconstructions in $\Omega_1$. First column: searched conductivity. Second column: reconstructed conductivity with unperturbed data ($\delta = 0$). Third column: reconstructed conductivity with perturbed data ($\delta = 2.0\%$).

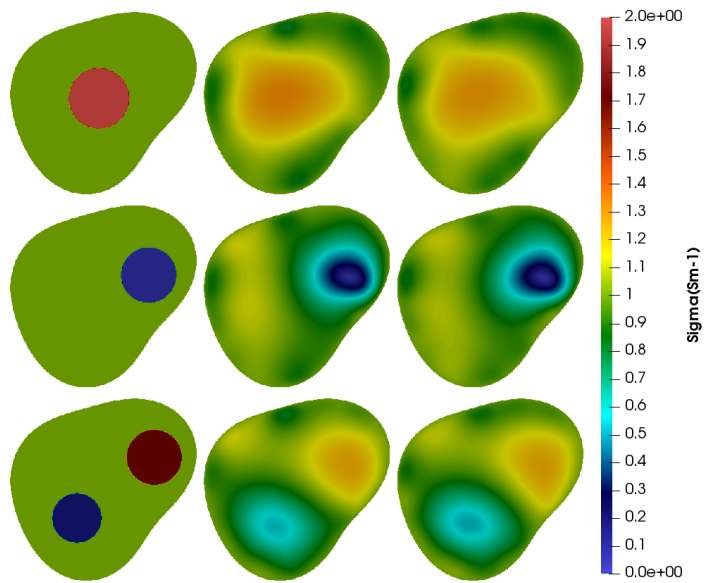

FIGURE 9. Reconstructions in $\Omega_3$. First column: searched conductivity. Second column: reconstructed conductivity with unperturbed data ($\delta = 0$). Third column: reconstructed conductivity with perturbed data ($\delta = 2.0\%$).

The partial derivatives of the functional with respect to parameters $\Theta_k^1$ or $\Theta_k^2$ are easily derived from the shape derivative obtained in [15], and in particular the sampling formula [15, Corollary 3.4]:

**Lemma 4.6.** *Let* $(\tilde{u}, \tilde{U})$ *be the solution of* (2.5) *for some electrode pattern* $\tilde{I} \in \mathbb{R}_\diamond^M$, *and* $f = g = 0$. *For $k$ in $\{1, 2\}$ and $m$ in $\{1, \ldots, M\}$, let $\mathbf{x}_m^k = r(\Theta_m^k)[\cos(\Theta_m^k), \sin(\Theta_m^k)]$.*

*We have for $m \in \{1, \ldots, M\}$:*

$$(4.7) \quad \frac{\partial U}{\partial \Theta_m^1} \cdot \tilde{I} = \rho(\Theta_m^1)(U_m - u(\mathbf{x}_m^1))(\tilde{U}_m - \tilde{u}(\mathbf{x}_m^1)) \quad \text{and} \quad \frac{\partial U}{\partial \Theta_m^2} \cdot \tilde{I} = \rho(\Theta_m^2)(U_m - u(\mathbf{x}_m^2))(\tilde{U}_m - \tilde{u}(\mathbf{x}_m^2)),$$

We refer to Section 4.1 for the definition of $r$ and $\rho$, and postpone the proof of Lemma 4.6 to the Appendix. The main difficulty with this sampling formula are the terms $u(\mathbf{x}_m^k)$ and $\tilde{u}(\mathbf{x}_m^k)$, which corresponds to the value of the interior potential at the end-points of the electrodes, and are difficult to compute accurately for two reasons: first, from a theoretical point of view, this potential is continuous at the end points of the electrodes, but no more than that, and may varies rapidly which makes it difficult to compute accurately. Furthermore, as in our immersed method the electrodes are not exactly discretized, such value is not directly available. We approximate it by linear interpolations from values of the nearest interface points located inside the electrode. Nevertheless, we obtain sufficiently precise results to reconstruct the position of the electrodes from the measurements by minimizing the functional $F$ using a gradient descent algorithm.

We first present reconstructions with $\Omega$ a disc of radius 0.5 covered with 4 electrodes. We either perturbed the first angle or the second angle of the first electrode, the other electrodes being unperturbed, but we nevertheless let all the angles be free in our algorithm. As expected, the electrodes $E_2$, $E_3$ and $E_4$, being already at their correct positions, are almost not moved by the algorithm, while the position of the first electrode is retrieved as shown in the following table:

| Searched angles | | Starting angles | | Reconstructions | |
|---|---|---|---|---|---|
| $\Theta_1^1$ | $\Theta_1^2$ | $\Theta_1^1$ | $\Theta_1^2$ | $\Theta_1^1$ | $\Theta_1^2$ |
| $-2.35619$ | $-1.85619$ | $-3.141592$ | $-1.85619$ | $-2.356028$ | $-1.85619$ |
| $-3.141592$ | $-2.64$ | $-3.141592$ | $-1.85619$ | $-3.14159$ | $-2.636460$ |

In our final experiments, we suppose that the length of the electrodes is fixed and known. In that case, $\Theta^1$ and $\Theta^2$ are not independent, as they are linked by relation (4.1). The method is easily adapted to this situation: we simply decide that $\Theta_2$ is now a function of $\Theta^1$ given by relation (4.1), and $\mathscr{U}$ depends only on $\mathscr{I}$ and $\Theta^1$. The derivatives are then obtain by simple combination of Lemma 4.6, the chain rule, and the fact that equation (4.1) immediately implies

$$\frac{\partial \Theta_k^2}{\partial \Theta_l^1} = \delta_{kl} \frac{\rho(\Theta_l^1)}{\rho(\Theta_k^2)}.$$

For our final tests, we choose $\Omega_e = [-1, 1] \times [-1, 1]$, discretized by a $400 \times 400$ Cartesian grid, $N = 4$ and $\alpha = [0.8, 0.02, 0.001, 0.05, 0.001, 0.04, 0.001]$. The obtained results are displayed in the following tables:

| | Searched positions | | Start positions | | Reconstructed positions | |
|---|---|---|---|---|---|---|
| | $\Theta^1$ | $\Theta^2$ | $\Theta^1$ | $\Theta^2$ | $\Theta^1$ | $\Theta^2$ |
| $E_1$ | $-2.51327$ | $-1.93817$ | $-3.14159$ | $-2.51358$ | $-2.55152$ | $-1.97563$ |
| $E_2$ | $-0.94247$ | $-0.30516$ | $-1.57079$ | $-0.91229$ | $-0.90298$ | $-0.27141$ |
| $E_3$ | $0.62831$ | $1.24343$ | $0$ | $0.57253$ | $0.54018$ | $1.14782$ |
| $E_4$ | $2.19911$ | $2.84837$ | $1.57079$ | $2.18999$ | $2.08957$ | $2.72868$ |

| | Searched positions | | Starting positions | | Reconstructed positions | |
|---|---|---|---|---|---|---|
| | $\Theta^1$ | $\Theta^2$ | $\Theta^1$ | $\Theta^2$ | $\Theta^1$ | $\Theta^2$ |
| $E_1$ | $-2.35619$ | $-1.78157$ | $-3.14159$ | $-2.51358$ | $-2.24361$ | $-1.66646$ |
| $E_2$ | $-0.78539$ | $-0.17053$ | $-1.57079$ | $-0.91229$ | $-0.70862$ | $-0.10368$ |
| $E_3$ | $0.78539$ | $1.41055$ | $0$ | $0.57253$ | $0.83745$ | $1.46474$ |
| $E_4$ | $2.19911$ | $2.84837$ | $1.57079$ | $2.18999$ | $2.23475$ | $2.88755$ |

Appendix A.

**Conductivity reconstruction – Proof of Lemma 4.4.** We prove Lemma 4.4 with $P = 1$, that is with a single current pattern $I$. The extension to several current patterns is straightforward. We therefore fix $I \in \mathbb{R}_\diamond^M$ and $U_{meas} \in \mathbb{R}^M$ such that $U_{meas,1} = 0$.

Let $(u(\sigma, I), U(\sigma, I))$ be the solution of (2.5) for a given conductivity $\sigma \in L^\infty \cap H^1(\Omega)$ satisfying $\sigma \geq c > 0$ a.e. in $\Omega$. It is known that the mapping

$$M : \sigma \longrightarrow (u(\sigma, I), U(\sigma, I)),$$

is Frechet differentiable [32]. More precisely, for all $\delta\sigma$ in $L^\infty \cap H^1(\Omega)$ such that $\sigma + \delta\sigma \geq \tilde{c} > 0$ a.e. in $\Omega$, one has

$$M(\sigma + \delta\sigma) = M(\sigma) + (\delta u, \delta U) + o(\|\delta\sigma\|_{L^\infty(\Omega)}),$$

where $\delta u \in H^1(\Omega)$ and $\delta U \in \mathbb{R}^M$, $\delta U_1 = 0$, are the only solution of the following variational problem: find $(u, U) \in H^1(\Omega) \times \mathbb{R}^M$ such that for all $(v, V) \in H^1(\Omega) \times \mathbb{R}^M$, one has

$$(A.1) \quad \int_\Omega \sigma \nabla u . \nabla v \, dx + \sum_m \int_{E_m} \xi_m (u - U_m)(v - V_m) \, ds(x) + \varepsilon U_1 V_1 = -\int_\Omega \delta\sigma \nabla u(\sigma, I) . \nabla v \, dx.$$

Equivalently, $(\delta u, \delta U)$ satisfies

$$\begin{cases} -\nabla \cdot (\sigma \nabla \delta u) = -\nabla \cdot (\delta\sigma \nabla u(\sigma, I)) & \text{in } \Omega \\ \sigma \partial_\nu (\delta u) = \xi_m (\delta U_m - \delta u) - (\delta\sigma) \partial_\nu u(\sigma, I) & \text{on } E_m \\ \sigma \partial_\nu (\delta u) = -(\delta\sigma) \partial_\nu u(\sigma, I) & \text{on } \partial\Omega \setminus \overline{E}, \\ \int_{E_m} \sigma \partial_\nu (\delta u) \, ds(x) + \varepsilon \delta U_1 = -\int_{E_m} (\delta\sigma) \partial_\nu (u(\sigma, I)) \, ds(x) & \text{for } m = 1, \ldots, M. \end{cases}$$

As a consequence, we have

$$(A.2) \quad F(\sigma + \delta\sigma) = F(\sigma) + (\delta U, U(\sigma, I_{input}) - U_{meas})_{\mathbb{R}^M} + \epsilon(\delta\sigma, \sigma - \sigma_*)_{H^1(\Omega)} + o(\|\delta\sigma\|_{L^\infty \cap H^1(\Omega)}).$$

We now define $(w, W)$ as the unique solution of (2.5) with $f = g = 0$ and input current $I = U(\sigma, I_{input}) - U_{meas}$, that is

$$w = u(\sigma, U(\sigma, I_{input}) - U_{meas}), \quad W = U(\sigma, U(\sigma, I_{input}) - U_{meas}).$$

Using (A.1) and the variational formulation of problem (2.5), we obtain that

$$(\delta U, U(\sigma, I_{input}) - U_{meas})_{\mathbb{R}^M} = -\int_\Omega \delta\sigma \nabla u(\sigma, I) \cdot \nabla w \, dx.$$

As a consequence, for $\delta\sigma$ in $L^\infty \cap H^1(\Omega)$ such that $\sigma + \delta\sigma \geq \tilde{c} > 0$ in $\Omega$, and $t \in (0, 1)$, we have

$$F(\sigma + \delta\sigma) = F(\sigma) + t \int_\Omega \left( \delta\sigma(\epsilon(\sigma - \sigma_*) - \nabla u(\sigma, I) . \nabla w) + \epsilon \nabla(\delta\sigma) . \nabla(\sigma - \sigma_*) \right) dx + o\left(t \|\delta\sigma\|_{L^\infty \cap H^1(\Omega)}\right)$$

Our aim is now to choose $\delta\sigma$ such that $F(\sigma + t\delta\sigma) \leq F(\sigma)$ for $t > 0$ small enough. The previous formula shows that it is sufficient to choose $\delta\sigma$ such that

$$\int_\Omega \left( \delta\sigma(\epsilon(\sigma - \sigma_*) - \nabla u(\sigma, I) . \nabla w) + \epsilon \nabla(\delta\sigma) . \nabla(\sigma - \sigma_*) \right) dx < 0.$$

Let us define $v$ as the unique function of $H_0^1(\Omega)$ such that for all $\tilde{v}$ in $H_0^1(\Omega)$,

$$(A.3) \quad \int_\Omega (\nabla v \cdot \nabla \tilde{v} + v \tilde{v}) \, dx = -\int_\Omega (\mathbf{G} \cdot \nabla \tilde{v} + f \tilde{v}) \, dx,$$

where

$$\mathbf{G} = \epsilon \nabla(\sigma - \sigma_*), \quad f = \epsilon(\sigma - \sigma_*) - \nabla u(\sigma, I) . \nabla w.$$

With this definition, we have

$$\int_\Omega \left( v(\epsilon(\sigma - \sigma_*) - \nabla u(\sigma, I) . \nabla w) + \epsilon \nabla v \cdot \nabla(\sigma - \sigma_*) \right) dx = -\|v\|_{H^1(\Omega)}^2 < 0.$$

Therefore, assuming that $v$ belongs to $L^\infty(\Omega)$, we can choose $\delta\sigma = v$, which ends the proof of Lemma 4.4.

*Remark* A.1. Clearly, $v$ solution of (A.3) might fail to be in $L^\infty(\Omega)$. However, from a numerical point of view, we will obtain an approximation of $v$ which will, by construction, be in $L^\infty(\Omega)$. Such approximation turns out to be a good direction descent, as shown in the presented reconstructions.

**Electrode position reconstruction – Proof of Lemma 4.6.** Lemma 4.6 is a consequence of a more general result, that is the Frechet derivability of the measurement map with respect to small perturbation of the boundary of $\Omega$, obtained in [15]. We start by very briefly recalling the results proved in this paper.

For $h \in C^1(\partial\Omega, \mathbb{R}^d)$, we set $F[h] : x \in \partial\Omega \mapsto x + h(x)$ and

$$\partial\Omega_h = F[h](\partial\Omega).$$

For $h$ small enough, $\partial\Omega_h$ is the boundary of a smooth domain $\Omega_h$, which of course is a perturbation of $\Omega$. For $h$ still small enough, $\Omega_h$ is covered by $M$ well separated electrodes $E_{m,h}$ defined by

$$E_{m,h} = \{x + h(x),\ x \in E_m\}.$$

As a consequence, for any $I \in \mathbb{R}_\diamond^M$, we can denote $(u(h), U(h))$ the solution of (2.5) with $(\Omega, E_m)$ replaced by $(\Omega_h, E_{m,h})$, and $f = g = 0$, and define the measurement map

$$R : (h, I) \in \mathscr{B}_d \times \mathbb{R}_\diamond^M \mapsto U(h),$$

with

$$\mathscr{B}_d = \{h \in C^1(\partial\Omega, \mathbb{R}^n); \|h\|_{C^1(\partial\Omega, \mathbb{R}^d)} < d\},$$

$d$ being a small enough fixed constant. The main result of [15] is the Frechet derivability of $R$ with respect to $h$ ([15, Theorem 3.2]). A sampling formula is also provided ([15, Corollary 3.4]), allowing to easily compute this derivative. These results are gathered in the following Theorem.

**Theorem A.2.** *The operator $R$ is Fréchet differentiable at the origin with respect to the first variable. In other words, there exists a bounded bilinear operator*

$$R' : C^1(\partial\Omega; \mathbb{R}^d) \times \mathbb{R}_\diamond^M \mapsto \mathbb{R}^M$$

*such that*

$$\lim_{h \to 0} \frac{1}{\|h\|_{C^1}} \|R(h, \cdot) - R(0, \cdot) - R'h\|_\mathscr{L} = 0.$$

*Furthermore, one can compute $R'h$ component by component using the following sampling formula: for $I$ and $\tilde{I}$ in $\mathbb{R}_\diamond^M$, if one denotes $(u, U)$ and $(\tilde{u}, \tilde{U})$ the corresponding solutions to* (2.5)*, then*

$$\left((R'h)I\right) \cdot \tilde{I} = -\sum_{m=1}^M \int_{\partial E_m} (h \cdot \nu_{\partial E_m})(U_m - u)(\tilde{U}_m - \tilde{u})ds(x)$$

$$- \sum_{m=1}^M \frac{1}{z_m} \int_{E_m} h_\nu \left((d-1)(U_m - u)H - \partial_\nu u\right)(\tilde{U}_m - \tilde{u})ds(x) - \int_{\partial\Omega} h_\nu (\sigma\nabla u)_\tau (\nabla\tilde{u})_\tau ds(x).$$

*Here, $H$ is the mean curvature, $\nu$ (resp. $\nu_{\partial E_m}$) is the outward normal vector to $\Omega$ (resp. to $E_m$), and $h_\nu = h \cdot \nu$.*

Now that we have all the necessary tools, we can prove Lemma 4.6. In the following, we focus on $\Theta^2$, the computations with $\Theta^1$ being similar. We also do the computation for a single input current $I \in \mathbb{R}_\diamond^M$, the extension to $P$ input currents being straightforward. Finally, for any $x \in \partial\Omega$, we define $\theta(x)$ as the unique element of $[-\pi, \pi[$ such that $x = r(\theta(x))\mathbf{u}(\theta(x))$.

Let $m \in \{1, \ldots, M\}$ be fixed, and $\eta \in \mathbb{R}$ sufficiently small such that

$$\Theta_m^1 < \Theta_m^2 + \eta < \Theta_{m+1}^1,$$

where, to simplify notations, we set $\Theta^1_{M+1} = \Theta^1_1 + 2\pi$. We define $\Theta^2_\eta$ by
$$\Theta^2_{\eta,k} = \Theta^2_k + \delta_{km}\,\eta, \ \forall k \in \{1, \ldots, M\}.$$
We define, for a fixed constant $c > 0$,
$$\mathscr{H}_\eta = \left\{ h \in C^1(\partial\Omega; \mathbb{R}^2), \ \left| \begin{array}{ll} h(\Theta^k_1) = 0 & \forall k \in \{1, \ldots, M\}, \\ h(\Theta^k_2) = \delta_{km}\eta & \forall k \in \{1, \ldots, M\}, \\ \|h\|_{C^1} \leq c\eta. & \end{array} \right. \right\}$$
The set $\mathscr{H}_\eta$ is non empty: indeed, take $\chi \in C_c^\infty(\mathbb{R}; [0,1])$ such that $\mathrm{supp}(\chi) \in \,]-1,1[$ and $\chi(0) = 1$ and define
$$\varphi_\eta(\theta) = \eta\chi\Big(\frac{\theta - \Theta^2_m}{\Delta}\Big)$$
where
$$\Delta = \max\left(\frac{\Theta^2_m - \Theta^1_m}{2}, \frac{\Theta^1_{m+1} - \Theta^2_m}{2}\right).$$
Then set $\omega_\eta : \mathbb{R} \mapsto \mathbb{R}$ $2\pi$-periodic, defined by
$$\omega_\eta : \theta \in [\Theta^1_1, \Theta^1_1 + 2\pi[ \longrightarrow \theta + \varphi_\eta(\theta).$$
The, if the constant $c$ is chosen sufficiently large, the vector field
$$h : x \in \partial\Omega \longrightarrow r(\omega_\eta(\theta(x)))\mathbf{u}(\omega_\eta(\theta(x))) - r(\theta)\mathbf{u}(\theta(x)),$$
belongs to $\mathscr{H}_\eta$.

As, for $\eta$ small enough, for $h_\eta$ in $\mathscr{H}_\eta$, one has $\partial\Omega_{h_\eta} = \partial\Omega$, $E_{k,h_\eta} = E_k$ for all $k \in \{1, \ldots, M\}$, $k \neq m$, and
$$E_{m,h_\eta} = \left\{ r(\theta)\mathbf{u}(\theta), \ \theta \in [\Theta^1_m, \Theta^2_m + \eta] \right\},$$
we have $U(\Theta^1, \Theta^2_\eta, I) = R(h_\eta, I)$. Now, cumbersome computations show that, for all $h_\eta$ in $\mathscr{H}_\eta$, for all $x$ in $\partial\Omega$, one has
$$h_\eta(x) = \Big(\theta(x + h_\eta(x)) - \theta(x)\Big)\rho(\theta(x))\tau(x) + o(\|h_\eta\|_{C^1}).$$
Applying Theorem A.2 with $d = 2$, immediately gives that for all $\tilde{I} \in \mathbb{R}^M_\diamond$ and $(\tilde{u}, \tilde{U})$ corresponding solution of (2.5), one has
$$(R'(h_\eta)I) \cdot \tilde{I} = \eta\,\rho(\Theta^2_m)(U_m - u(\mathbf{x}^2_m))(\tilde{U}_m - \tilde{u}(\mathbf{x}^2_m)) + o(\|h_\eta\|_{C^1}).$$
As
$$\Big(U(\Theta^1, \Theta^2_\eta, I) - U(\Theta^1, \Theta^2, I)\Big) \cdot \tilde{I} - \eta\,\rho(\Theta^2_m)(U_m - u(\mathbf{x}^2_m))(\tilde{U}_m - \tilde{u}(\mathbf{x}^2_m))$$
$$= \Big(R(h_\eta, I) - R(0, I) - R'(h_\eta)I\Big) \cdot \tilde{I} + o(\|h_\eta\|_{C^1}),$$
Theorem A.2 directly implies Lemma 4.6.

**Acknowledgements.** This work has been partially supported by the ANR LabEx CIMI (under grant ANR-11-LABX-0040) within the French State Programme "Investissements d'Avenir".

Experiments presented in this paper were carried out using the PlaFRIM experimental testbed, supported by Inria, CNRS (LABRI and IMB), Université de Bordeaux, Bordeaux INP and Conseil Régional d'Aquitaine.

(J. Dardé) Institut de Mathématiques de Toulouse, UMR 5219, Université de Toulouse, CNRS, UPS, 118 route de Narbonne, 31062 Toulouse Cédex 9, France.

(N. Nasr) Univ. Bordeaux, CNRS, INRIA, Bordeaux INP, IHU-LIRYC, IMB, UMR 5251, F-33400 Talence, France

(L. Weynans) Univ. Bordeaux, CNRS, INRIA, Bordeaux INP, IHU-LIRYC, IMB, UMR 5251, F-33400 Talence, France

*Email address*, J. Dardé: `jeremi.darde@math.univ-toulouse.fr`

*Email address*, N. Nasr: `niami.nasr@u-bordeaux.fr`

*Email address*, L. Weynans: `lisl.weynans@math.u-bordeaux.fr`